\documentclass[11pt]{article}
\usepackage{color}
\usepackage[colorlinks=true]{hyperref}
\usepackage[margin=1in]{geometry}
\usepackage{amssymb}
\usepackage{amsmath}
\usepackage{amsthm}
\usepackage{mathtools}
\usepackage{epsfig, graphicx}
\usepackage{subfigure} 
\usepackage{ulem}

\newtheorem{theorem}{Theorem}[section]

\numberwithin{equation}{section}

\allowdisplaybreaks 
\begin{document}
\title{Tensor Neural Network Based Machine Learning Method for Elliptic Multiscale 
Problems\footnote{This work was supported  by  the Beijing Natural Science Foundation (No. Z200003), 
the National Key Laboratory of Computational Physics of China (6142A05230501), the National Center for 
Mathematics and Interdisciplinary Science, CAS.}}
\author{Zhongshuo Lin\footnote{LSEC, NCMIS, Institute
of Computational Mathematics, Academy of Mathematics and Systems
Science, Chinese Academy of Sciences, Beijing 100190,
China,  and School of Mathematical Sciences, University
of Chinese Academy of Sciences, Beijing 100049, China (linzhongshuo@lsec.cc.ac.cn).}, \ \
Haochen Liu\footnote{LSEC, NCMIS, Institute
of Computational Mathematics, Academy of Mathematics and Systems
Science, Chinese Academy of Sciences, Beijing 100190,
China,  and School of Mathematical Sciences, University
of Chinese Academy of Sciences, Beijing 100049, China (liuhaochen@lsec.cc.ac.cn).},
\ \  and \ \
Hehu Xie\footnote{LSEC, NCMIS, Institute
of Computational Mathematics, Academy of Mathematics and Systems
Science, Chinese Academy of Sciences, Beijing 100190,
China,  and School of Mathematical Sciences, University
of Chinese Academy of Sciences, Beijing 100049, China (hhxie@lsec.cc.ac.cn).}}

\date{}
\maketitle

\begin{abstract}
In this paper, we introduce a type of tensor neural network based machine learning method to solve 
elliptic multiscale problems. Based on the special structure, we can do the
direct and highly accurate high dimensional integrations 
for the tensor neural network functions without Monte Carlo process.
Here, with the help of homogenization techniques, 
the multiscale problem is first transformed to the high dimensional limit problem with
reasonable accuracy. Then, based on the tensor neural network, we design a type of
machine learning method to solve the derived high dimensional limit problem.
The proposed method in this paper brings a new way to design numerical methods for computing 
more general multiscale problems with high accuracy. 
Several numerical examples are also provided to validate the accuracy of the proposed numerical methods. 


\end{abstract}

\section{Introduction}
In modern science and engineers, many real-world problems involve multiple 
spatial and temporal scales, posing significant challenges for traditional 
numerical methods. Many processes, such as materials 
with varying microstructures or multiscale physical phenomena, 
require the consideration of multiple scales to capture their behavior accurately.  
Multiscale partial differential equations (MPDEs) arise when these scales can 
not be fully resolved, necessitating the development of efficient numerical techniques. 
Numerical methods for MPDEs aim to accurately capture the behavior of systems 
with disparate scales while maintaining computational efficiency. 
Traditional methods, such as finite difference or finite element methods, 
fail to capture the dynamics of multiscale phenomena due to limitations in 
grid resolution or excessive computational costs. 
To overcome these challenges, special algorithms have been developed that combine 
insights from both deterministic and stochastic methods.

One commonly used approach is the concept of homogenization, whereby the MPDE is 
approximated by a simpler, effective equation at the macroscopic scale. 
Homogenization strategies aim to extract the macroscopic 
behavior by averaging over the fine-scale features of the system. 
These methods are particularly useful when the fine-scale information is not crucial 
for the global behavior of the problem \cite{AllaireBriane,Cao, CaoCui, JikovKozlovOleinik}. 
Multiscale finite element methods (MsFEM) represent another class of 
numerical algorithms for MPDEs \cite{ChenHou, EfendievHou,HouWu}. 
These methods typically employ enriched finite element spaces that capture the behavior 
of the solution on various scales. By augmenting the standard basis functions with additional 
localized functions, MsFEM can accurately represent the behavior of the solution across multiple scales.
Heterogeneous multiscale method (HMM) is a powerful computational approach designed to tackle problems that 
exhibit heterogeneity across multiple scales \cite{HMM_Book,EMingZhang}. By integrating information from 
fine-scale simulations and coarse-scale models, 
it strikes a balance between accuracy and computational efficiency. Gamblets \cite{Owhadi,OwhadiZhang} 
(operator-adapted wavelets satisfying three desirable properties: scale orthogonality, 
well conditioned multi-resolution decomposition, and localization) provide 
a natural multiresolution decomposition ensuring robustness for multiscale problems. 
Furthermore, there have developed fast gamblet solvers 
to achieve grid-size accuracy for elliptic problems. Recently, a type of tensor finite element 
method is designed for solving the 
MPDEs \cite{HarbrechtSchwab,HoangSchwab} which brings another way to deal with the 
MPDEs with scale separation. 
Furthermore, stochastic methods, such as Monte Carlo methods or stochastic finite element methods \cite{AbdulleBarthSchwab}, 
have been successfully employed to tackle MPDEs 
characterized by uncertain or random inputs. These methods incorporate random fluctuations 
in the system's parameters or initial/boundary conditions, providing a statistical 
representation of the solution that encompasses various scales.

Researchers also studied the data-driven approaches to solve multiscale problems.
Recently, there appear several methods such as Fourier neural operator (FNO) \cite{FNO}, 
Galerkin transformer (GT) \cite{GT}, deep operator network (DeepONet) \cite{DeepONet},  
which are designed to directly learn the operator (mapping) between infinite dimensional 
spaces for partial differential equations (PDEs) 
by taking advantages of the enhanced expressibility of deep neural 
networks and advanced architectures such as feature embedding, channel mixing and self-attentions.  
However, for multiscale problems, most existing operator learning schemes essentially capture the 
smooth part of the solution space, and how to resolve the intrinsic multiscale features 
remains to be a major challenge.
In \cite{LiuXuZhang},  a type of hierarchical transformer neural network (HT-Net) with 
scale adaptive interaction
range, such that the features can be computed in a nested manner and with a controllable linear
cost. Self-attentions over a hierarchy of levels can be used to encode and decode the multiscale
solution space over all scale ranges. 
Physics-informed neural network (PINN) is a neural network approach to 
solve partial differential equation 
\cite{RaissiPerdikarisKarniadakis}.   
Compared with the classical numerical methods, PINN is a mesh free method and can 
interpret the solution without a mesh.  
There also exist several PINN based numerical methods for solving multiscale problems in 
\cite{ChungEfendierLeungPunZhang, ChungLeungPunZhang, LeungLinZhang, 
WangWangPerdikaris, ZhangChungEfendierLeung}.

So far, the machine learning based method for MPDEs are designed based on 
sampling points or Monte Carlo process.  
Unfortunately, the sampling step destroy the accuracy of the  machine learning method. 
Recently, there appears a type of tensor neural network (TNN) and the corresponding 
machine learning method which are designed to solve high dimensional 
problems with high accuracy \cite{WangJinXie,WangLiaoXie,WangXie}. 
The most important property  of TNN is that the corresponding high dimensional functions 
can be easily integrated with high accuracy and high efficiency. Then the deduced machine 
learning method can achieve high accuracy for solving high dimensional problems. 
The reason is that the integration of TNN functions can be separated into one dimensional 
integrations which can be done by the classical quadrature schemes with high accuracy. 
The tensor structure here is built by using the correlated one-dimensional-input 
feedforward neural network and tensor product way to construct the trial function space. 
With the powerful adaptiveness of neural network based machine learning methods, 
the trial functions possess larger expressive power and easier implementation than that of the classical 
numerical methods which are built with the fixed basis.

As a new application, the aim of this paper is to propose a type of TNN based machine learning method to 
solve MPDEs with scale separation. 
First, the MPDEs are expanded with the multiscale series which is determined by 
some high dimensional partial differential equations (PDEs) with one scale. 
These high dimensional PDEs are called limit problems.  
Then, TNN based machine learning method is adopted to solve these 
high dimensional limit problems to obtain the solution of MPDEs.  
Compared with existing machine learning methods, the proposed method here 
can reach high accuracy as the classical numerical methods for solving low dimensional problems.   
The essential reason is that the high dimensional integrations of TNN functions 
can be computed with high accuracy and high efficiency, which also makes it easy 
to guarantee the accuracy for the constraint conditions in the problems.

An outline of the paper goes as follows. In Section \ref{Section_MPDEs}, 
we review the homogenization method for the concerned multiscale PDEs. 
We then introduce the structure of the TNN and its approximation property 
in Section \ref{Section_TNN}. Section \ref{homo-TNN} is devoted to proposing the 
TNN-based machine learning method for solving MPDEs based on 
the high dimensional limit problems. Especially, the numerical integration 
method for the high dimensional integration in the loss functions is 
designed in Section \ref{Section_Integration}.
Some numerical examples are provided in Section \ref{Section_Numerical}
to show the validity and efficiency of the proposed numerical methods in this paper.
Some concluding remarks are given in the last section.

\section{Solving multiscale partial differential equations}\label{Section_MPDEs}
In this section, we introduce the well-known homogenization technique for solving 
the multiscale PDEs. We expand the MPDEs into multiscale series 
which are determined by some high dimensional PDEs. Then these high dimensional PDEs 
will be solved by the TNN based the TNN based machine learning method in the next sections.

As an example, we consider the following elliptic multiscale problem: 
Find $u^\varepsilon(\mathbf{x})\in H_0^1(\Omega)$ such that
\begin{eqnarray}\label{problem}
\left\{
\begin{array}{rcl}
-{\rm div}\left(A^\varepsilon \nabla u^\varepsilon\right)&=& f,\ \ \ {\rm in}\ \Omega,\\
u^\varepsilon &=&0,\ \ \ {\rm on}\ \partial\Omega.
\end{array}
\right.
\end{eqnarray}
The $d\times d$ matrix $A^\varepsilon = (A_{ij})_{1\leq i,j\leq d}$ 
is assumed to depend on the $\varepsilon$ with multiscales 
in the following sense: there exist $K$ positive scale functions 
$\varepsilon_1, \cdots, \varepsilon_{K}$ of $\varepsilon$
such that
\begin{eqnarray}\label{Scale_Separation_Condition}
\lim_{\varepsilon\rightarrow 0}\frac{\varepsilon_{k+1}}{\varepsilon_k} = 0,
\ \  {\rm for}\ k = 1, \cdots, K-1,
\end{eqnarray}
and for all $\mathbf x\in \Omega$ and all $0<\varepsilon_k<1$
$$
A^\varepsilon(\mathbf x) = A\left(\mathbf x, \frac{\mathbf x}{\varepsilon_1}, 
\cdots, \frac{\mathbf x}{\varepsilon_{K}}\right).
$$
In this paper, the matrix $A$ is assumed to be bounded and uniformly positive definite, i.e., 
there is a constant $\gamma>0$ such that 
\begin{eqnarray}\label{Elliptic_Assumption}
\gamma|\xi|^2 \leqslant \xi^{\top} A\left(\mathbf{x}, \mathbf{y}_1, \ldots, 
\mathbf{y}_{K}\right) \xi \leqslant \gamma^{-1}|\xi|^2,\ \ \ \forall \xi \in \mathbb{R}^{d}, 
\end{eqnarray}
for all $\mathbf{x} \in \Omega$ and $\mathbf{y}_k := 
\frac{\mathbf{x}}{\varepsilon_k} \in Y_k = [0,1]^d$, $k=1, \cdots, K$. 
Furthermore, we assume that $A$ is periodic with respect to 
$\mathbf{y}_k$ in domain $Y_k = [0,1]^d$. 
In order to deduce the limit problems, we seek the following asymptotic expansion 
for $u^{\varepsilon}(\mathbf{x})$:
\begin{eqnarray}\label{u_e_expansion}
u^{\varepsilon}(\mathbf{x}) = u_0(\mathbf{x}) + \varepsilon_1 u_1(\mathbf{x}, \mathbf{y}_1) 
+ \varepsilon_2 u_2(\mathbf{x}, \mathbf{y}_1, \mathbf{y}_2) + \cdots 
+ \varepsilon_{K}u_{K}(\mathbf{x}, \mathbf{y}_1, \cdots, \mathbf{y}_{K}),
\end{eqnarray}
where $u_k(\mathbf{x}, \mathbf{y}_1, \cdots,\mathbf{y}_k), k=1, \cdots, K$, are periodic in $\mathbf{y}_j$, 
$j = 1,\cdots,k, k\geq 1$. In order to guarantee the uniqueness of the asymptotic expansion,  
the functions $u_k(\mathbf{x}, \mathbf{y}_1, \cdots,\mathbf{y}_k), k =1, \cdots, K$, are assumed to 
satisfy the following conditions: For $k = 1,\cdots, K, $
\begin{eqnarray}\label{unique_cond}
\int_{{Y}_k} u_k(\mathbf{x}, \mathbf{y}_1, \cdots,\mathbf{y}_k) d\mathbf{y}_k =0, \ \ 
 \forall  \mathbf{x}\in\Omega,\ \ \forall \mathbf{y}_j\in Y_j, \  j< k.
\end{eqnarray}
For the sake of notation, in this paper, we set $\mathbf{y}_0:=\mathbf{x}$ 
and $\varepsilon_k := \varepsilon^k$, where $0<\varepsilon<1$. 
Then it is easy to verify that:
\begin{eqnarray}\label{L}
\text{div}(A(\mathbf{x}, \mathbf{y}_1, \cdots, \mathbf{y}_K)\nabla u^{\varepsilon}(\mathbf{x})) = \sum_{i,j,k=0,\cdots,K}\varepsilon^{k-i-j}\nabla_{\mathbf{y}_i}\cdot
\left(A\nabla_{\mathbf{y}_j}u_k\right).
\end{eqnarray}
Notice that $-{\rm div}(A(\mathbf{x}, \mathbf{y}_1, \cdots, 
\mathbf{y}_K)\nabla u^{\varepsilon}(\mathbf{x})) = f(\mathbf{x})$. 
With the help of (\ref{L}), equating the terms with the same power of $\varepsilon$ 
on both sides of (\ref{problem}) leads to the following equations 
\begin{eqnarray}\label{terms_equal}
\left\{
\begin{array}{rcl}
\sum\limits_{\stackrel{i,j,k=0,\cdots,K}{k-i-j=-n}}\nabla_{\mathbf{y}_i}\cdot
\left(A\nabla_{\mathbf{y}_j}u_k\right)& = & 0,  \ \ \forall n = 1, \cdots, K,\\
&&\\
\sum\limits_{\stackrel{i,j,k=0,\cdots,K}{k-i-j=0}}\nabla_{\mathbf{y}_i}\cdot
\left(A\nabla_{\mathbf{y}_j}u_k\right)&=& -f, 
\end{array}
\right.
\end{eqnarray}
where we neglect all the $\mathcal{O}(\varepsilon)$ terms. 
Using the expansion (\ref{u_e_expansion}) and the fact that $A$, 
as well as $u_k(\mathbf{x}, \mathbf{y}_1, \cdots,\mathbf{y}_k)$, 
is periodic with respect to $\mathbf{y}_k, k=1, \cdots,K$, we have
\begin{eqnarray}
\nabla_{\mathbf{y}_j}u_k=0,&&\ \  \forall j>k, \label{Equality_1}\\
\int_{Y_i} \nabla_{\mathbf{y}_i}\cdot (A\nabla_{\mathbf{y}_j}u_k) d\mathbf{y}_i = 0,
&&\ \ \forall i>0,j \geq 0,k\geq 0. \label{Equality_2}
\end{eqnarray}
Based on (\ref{Equality_1})-(\ref{Equality_2}), (\ref{terms_equal}) 
can be simplified by integrating both sides to the following equations:
\begin{eqnarray}
\nabla_{\mathbf{y}_K}\cdot\left(A\left(\sum_{i=0,\cdots,K}\nabla_{\mathbf{y}_i}u_i\right)\right)&=& 0, \label{terms_equal_simple1}\\
\int_{{Y}_K}\nabla_{\mathbf{y}_{K-1}}\cdot
\left(A\left(\sum_{i=0,\cdots,K}\nabla_{\mathbf{y}_i}u_i\right)\right)d\mathbf{y}_K&=& 0, \label{terms_equal_simple2}\\
\int_{{Y}_{K}}\int_{{Y}_{K-1}}\nabla_{\mathbf{y}_{K-2}}\cdot
\left(A\left(\sum_{i=0,\cdots,K}\nabla_{\mathbf{y}_i}u_i\right)\right)
d\mathbf{y}_{K-1}d\mathbf{y}_K&=& 0, \label{terms_equal_simple3}\\
\vdots&&  \nonumber\\
\int_{{Y}_{K}}\cdots\int_{{Y}_{1}}\nabla_{\mathbf{x}}\cdot
\left(A\left(\sum_{i=0,\cdots,K}\nabla_{\mathbf{y}_i}u_i\right)\right)
d\mathbf{y}_{1}\cdots d\mathbf{y}_K& = & -f.\label{terms_equal_simple4}
\end{eqnarray}
In equation (\ref{terms_equal_simple1}), if we view 
$m_j^K=(\mathop{\sum}_{i=0,\cdots,K-1}\nabla_{\mathbf{y}_i}u_i)_j$ 
as a constant coefficient, $u_K$ can be written as
\begin{eqnarray}\label{uK-uK-1}
u_K(\mathbf{x}, \mathbf{y}_1, \cdots, \mathbf{y}_K) = \boldsymbol{\chi}^K(\mathbf{x}, \mathbf{y}_1, 
\cdots, \mathbf{y}_K)\cdot M^K(\mathbf{x}, \mathbf{y}_1, \cdots, \mathbf{y}_{K-1}),
\end{eqnarray}
where $\boldsymbol{\chi}^K = [\mathbf{\chi}^K_1, \mathbf{\chi}^K_2, \cdots, \mathbf{\chi}^K_d]^\top$ 
and $ M^K = [m_1^K, m_2^K, \cdots, m_d^K]^\top$ with 
$\mathbf{\chi}^K_j$, $j=1,\cdots,d$, being defined as the solutions of the following cell problems:
\begin{eqnarray}\label{cell-prob}
\left\{
\begin{array}{lcl}
-\nabla_{\mathbf{y}_K}\cdot\Big(A\nabla_{\mathbf{y}_K}\mathbf{\chi}^K_j\Big)= \nabla_{\mathbf{y}_K}\cdot(A\mathbf{e}_j), 
\ \ j=1,\cdots,d,&&\\
\mathbf{\chi}^K_j \ \text{is periodic in}\   Y_K=[0,1]^d \ \text{with mean} \ 0,&&
\end{array}
\right.
\end{eqnarray}
where $(\mathbf{e}_j)_i = \delta_{ji}, j = 1,\cdots,d$. 
Then, (\ref{terms_equal_simple2})-(\ref{terms_equal_simple4}) can be written
as following form
\begin{eqnarray}
\nabla_{\mathbf{y}_{K-1}}\cdot\left(A^{*}_{K-1}\left(\sum_{i=0,\cdots,K-1}
\nabla_{\mathbf{y}_i}u_i\right)\right)&=& 0, \label{newform1}\\
\int_{{Y}_{K-1}}\nabla_{\mathbf{y}_{K-2}}\cdot\left(A^{*}_{K-1}
\left(\sum_{i=0,\cdots,K-1}\nabla_{\mathbf{y}_i}u_i\right)\right)d\mathbf{y}_{K-1}&=& 0, 
\label{newform2}\\
\vdots&&  \nonumber\\
\int_{{Y}_{K-1}}\cdots\int_{{Y}_{1}}\nabla_{\mathbf{x}}\cdot\left(A^{*}_{K-1}
\left(\sum_{i=0,\cdots,K-1}\nabla_{\mathbf{y}_i}u_i\right)\right)
d\mathbf{y}_{1}\cdots d\mathbf{y}_{K-1}& = & -f,\label{newform3}
\end{eqnarray}
where the homogenized matrix $A_{K-1}^* = \int_{{Y}_{K}} A(\mathbf{I}_d+\nabla_{\mathbf{y}_K}\boldsymbol{\chi}^K)d\mathbf{y}_{K}$. 
Here, $\mathbf{I}_d$ is the $d \times d$ identity matrix and 
$\nabla_{\mathbf{y}_K}\boldsymbol{\chi}^K$ is the transpose of the 
Jacobian matrix of $\boldsymbol{\chi}^K$ with respect to $\mathbf{y}_K$, i.e., 
$(\nabla_{\mathbf{y}_K}\boldsymbol{\chi}^K)_{ij} 
= \frac{\partial \chi^K_j}{\partial \mathbf{y}_{Ki}}, 1\leq i,j\leq d$.   
It is easy to observe that equations (\ref{newform1})-(\ref{newform3}) is a 
system similar to (\ref{terms_equal_simple1})-(\ref{terms_equal_simple4}) 
but with multiscale coefficient $A_{K-1}^*$ and one scale less. 
Therefore, we can iterate the process as above and arrive at the 
following homogenized equation for $u_0$:
\begin{eqnarray}\label{homo-homogenization}
-\nabla_{\mathbf{x}} \cdot \bigg( A^{*}_{0}\nabla_{\mathbf{x}}u_0 \bigg) & = & f.
\end{eqnarray}
Solving the homogenized equation (\ref{homo-homogenization}) for $u_0$, 
then following a back substitution process with
\begin{eqnarray}\label{uk-uk-1}
u_k(\mathbf{x}, \mathbf{y}_1, \cdots, \mathbf{y}_k) = \boldsymbol{\chi}^k(\mathbf{x}, 
\mathbf{y}_1, \cdots, \mathbf{y}_k)\cdot M^k(\mathbf{x}, \mathbf{y}_1, 
\cdots, \mathbf{y}_{k-1}), \quad k=1, \cdots, K,
\end{eqnarray}
we can obtain $u_1, \cdots, u_K$ in order, and therefore an approximation $u^{\varepsilon}(\mathbf{x}) 
= u_0(\mathbf{x}) + \varepsilon u_1(\mathbf{x}, \frac{\mathbf{x}}{\varepsilon}) 
+ \cdots + \varepsilon^K u_K(\mathbf{x}, \frac{\mathbf{x}}{\varepsilon}, \cdots, \frac{\mathbf{x}}{\varepsilon^K})$ for the 
problem (\ref{problem}).

In particular, for the case of $d=1$, the cell problem (\ref{cell-prob}) can be written as
\begin{eqnarray}\label{cell-1D}
-\frac{\partial}{\partial y_K}
\left(A\left(1+\frac{\partial\chi^K_1}{\partial y_K}\right)\right)=0.
\end{eqnarray}
Integrating on $Y_K,$ we can get
\begin{eqnarray}
A(x,y_1,\cdots,y_{K})\left(1+\frac{\partial\chi^K_1(x,y_1,\cdots,y_{K})}{\partial y_K}\right)=C(x,y_1,\cdots,y_{K-1}).
\end{eqnarray}
Dividing by $A(x,y_1,\cdots,y_{K})$ and integrating on $Y_K$ leads to following equation 
\begin{eqnarray}
C(x,y_1,\cdots,y_{K-1}) = \frac{1}{\int_{Y_K}\frac{1}{A(x,y_1,\cdots,y_{K})}dy_K}.
\end{eqnarray}
Hence, we obtain that
\begin{eqnarray}
A_{K-1}^* = \int_{{Y}_{K}} A\left(1+\frac{\partial\chi^K_1}{\partial y_K}\right)dy_{K} 
= C = \frac{1}{\int_{Y_K}\frac{1}{A}dy_K}.
\end{eqnarray}
The same technique can derive that
\begin{eqnarray}\label{homogenized-matrix-analytical}
A_{k-1}^* = \frac{1}{\int_{Y_k} \cdots \int_{Y_K} \frac{1}{A} dy_K \cdots dy_k}, 
\quad k = 1,2,\cdots, K
\end{eqnarray}
in the iteration process. In particular, $A_0^*$ has following form 
\begin{eqnarray}\label{homo-coefficient-analytical}
A_{0}^* = \frac{1}{\int_{Y_1} \cdots \int_{Y_K} \frac{1}{A} dy_K \cdots dy_1}.
\end{eqnarray}

In this paper, the TNN based machine learning method is developed 
for the particular case of $K=1, 2$. In the following, 
we take the case of $K = 1$ as an example to express the solving procedure.  
We consider the following elliptic multiscale problem in divergence form
\begin{eqnarray}\label{Multiscale_Problem_2}
-\operatorname{div}\left(A\left(\mathbf{x}, \frac{\mathbf{x}}{\varepsilon}\right) 
\nabla u^{\varepsilon}(\mathbf{x})\right)=f(\mathbf{x}).
\end{eqnarray}
The corresponding cell problems are: Find $\chi_j(\mathbf{x}, \mathbf{y}), j=1,\cdots,d$, 
such that 
\begin{eqnarray}
\left\{
\begin{array}{lcl}
-\nabla_{\mathbf{y}}\cdot\Big(A\nabla_{\mathbf{y}}\mathbf{\chi}_j\Big)
= \nabla_{\mathbf{y}}\cdot(A\mathbf{e}_j), \ \ j=1,\cdots,d,&&\\
\mathbf{\chi}_j \ \text{is periodic in}\   Y=[0,1]^d \ \text{with mean} \ 0.&&
\end{array}
\right.
\end{eqnarray}
Then the homogenized coefficient can be defined as follows 
\begin{eqnarray}
A^* = \int_{Y} A(\mathbf{I}_d +\nabla_{\mathbf{y}}\boldsymbol{\chi})d\mathbf{y}.
\end{eqnarray}
The homogenized equation for $u_0$ is: Find $u_0\in H_0^1(\Omega)$ such that 
\begin{eqnarray}
-\nabla_{\mathbf{x}}\cdot(A^*\nabla_{\mathbf{x}} u_0) = f.
\end{eqnarray}
Finally, we obtain $u_1$ by
\begin{eqnarray}
u_1(\mathbf{x},\mathbf{y}) = \boldsymbol{\chi}(\mathbf{x},\mathbf{y}) 
\cdot \nabla_{\mathbf{x}} u_0(\mathbf{x}),
\end{eqnarray}
and the multiscale approximation solution by 
\begin{eqnarray}
u^{\varepsilon}(\mathbf{x}) = u_0(\mathbf{x}) + 
\varepsilon u_1\left(\mathbf{x},\frac{\mathbf{x}}{\varepsilon}\right).
\end{eqnarray}

Proposition 3.3 of \cite{HoangSchwab} gives the convergence 
result of the homogenization method for the case of $K=1$. 
Assume that $A(\mathbf{x}, \mathbf{y}) \in C^{\infty}\left(\bar{\Omega}, 
C_{\text{per}}^{\infty}(Y)\right)_{\text{sym}}^{d \times d}$ and 
the homogenized solution $u_0(\mathbf{x})$ belongs to $H^2(\Omega)$. Then
\begin{eqnarray}\label{Convergence_L2}
\left\|u^{\varepsilon}-\left(u_0(\mathbf{x})+\varepsilon u_1\left(\mathbf{x}, 
\frac{\mathbf{x}}{\varepsilon}\right)\right)\right\|_{L^2(\Omega)} \leq C_1 \varepsilon,
\end{eqnarray}
\begin{eqnarray}\label{Convergence_H1}
\left\|u^{\varepsilon}-\left(u_0(\mathbf{x})+\varepsilon u_1\left(\mathbf{x}, 
\frac{\mathbf{x}}{\varepsilon}\right)\right)\right\|_{H^1(\Omega)} \leq C_2 \varepsilon^{1/2},
\end{eqnarray}
where the constant $C_1, C_2 > 0$ are independent of $\varepsilon$ 
but depend on $u_0$ and $u_1$.
Similar results are also presented in
\cite{BensoussanLionsPapanicolaou,JikovKozlovOleinik}.

\section{Tensor neural network architecture}\label{Section_TNN}
This section is devoted to introducing the TNN structure and some 
techniques to improve its stability. The approximation property and the computational complexity of 
related integration of TNN functions have been discussed and investigated in \cite{WangJinXie}.
In order to express clearly and facilitate the construction of the TNN method for solving
high dimensional PDEs, here, we will also elaborate on some important approximation properties.

\begin{figure}[htb!]
\centering
\includegraphics[width=16cm,height=12cm]{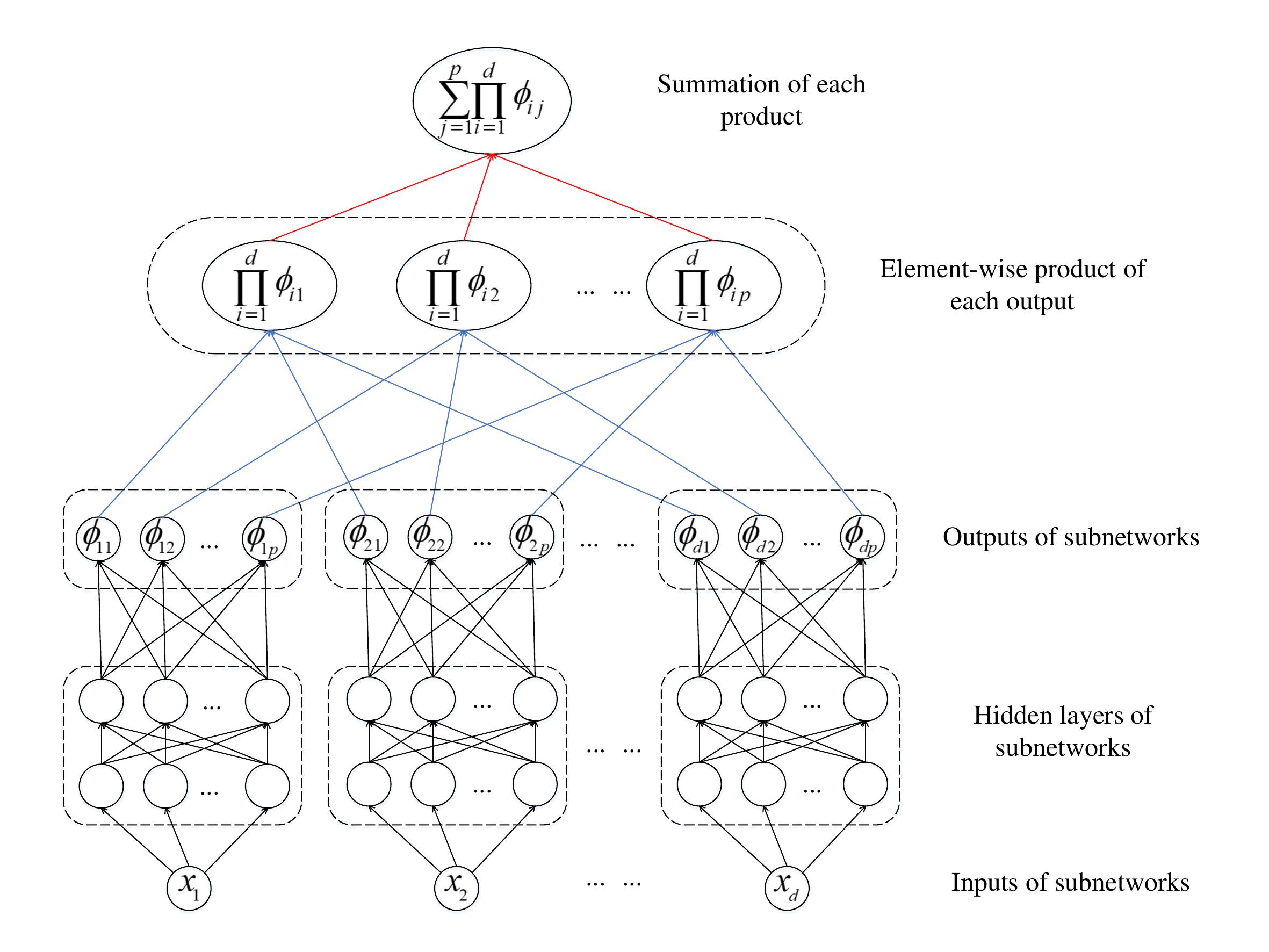}
\caption{Architecture of TNN. Black arrows mean linear transformation
(or affine transformation). Each ending node of blue arrows is obtained by taking the
scalar multiplication of all starting nodes of blue arrows that end in this ending node.
The final output of TNN is derived from the summation of all starting 
nodes of red arrows.}\label{TNNstructure}
\end{figure}

TNN is built by the tensor products of one dimension functions which come from
$d$ subnetworks with one-dimensional input and multidimensional output, where
$d$ denotes the spatial dimension of the concerned problems which will be solved by the
machine learning method in this paper. 
For each $i=1,2,\cdots,d$, we use $\Phi_i(x_i;\theta_i)=(\phi_{i,1}(x_i;\theta_i),
\phi_{i,2}(x_i;\theta_i), \cdots,\phi_{i,p}(x_i;\theta_i))$
to denote a subnetwork that maps a set $\Omega_i\subset\mathbb R$ to $\mathbb R^p$,
where $\Omega_i,i=1,\cdots,d,$ can be a bounded interval $(a_i,b_i)$, 
the whole line $(-\infty,+\infty)$ or the half line $(a_i,+\infty)$.
The number of layers and neurons in each layer, the selections of activation 
functions and other hyperparameters can be different in different subnetworks. 
TNN consists of $p$ correlated rank-one functions,
which are composed of the multiplication of $d$ one-dimensional input functions 
in different directions. Figure \ref{TNNstructure} shows the corresponding architecture of TNN.

In implementation, in order to improve the numerical stability, 
we normalize each $\phi_{i,j}(x_i)$
and use the following normalized-TNN structure:
\begin{eqnarray}\label{def_TNN_normed}
\Psi(x;c,\theta)=\sum_{j=1}^pc_j\widehat\phi_{1,j}(x_1;\theta_1)\widehat\phi_{2,j}(x_2;\theta_2)
\cdots\widehat\phi_{d,j}(x_d;\theta_d)
=\sum_{j=1}^pc_j\prod_{i=1}^d\widehat\phi_{i,j}(x_i;\theta_i),
\end{eqnarray}
where each $c_j$ is scaling parameter which describes the length of each rank-one 
function, $c=\{c_j\}_{j=1}^{p}$ is a set of trainable parameters, 
$\{c,\theta\}=\{c,\theta_1,\cdots,\theta_d\}$
denotes all parameters of the whole architecture.
For $i=1,\cdots,d,j=1,\cdots,p$, $\widehat\phi_{i,j}(x_i,\theta_i)$ is a $L^2$-normalized 
function as follows:
\begin{eqnarray}\label{eq_phi_normed}
\widehat\phi_{i,j}(x_i,\theta_i)
=\frac{\phi_{i,j}(x_i,\theta_i)}{\|\phi_{i,j}(x_i,\theta_i)\|_{L^2(\Omega_i)}}.
\end{eqnarray}
For simplicity of notation, $\phi_{i,j}(x_i,\theta_i)$ denotes the 
normalized function in the following parts.

The TNN architecture (\ref{def_TNN_normed}) and the architecture defined in \cite{WangJinXie} 
are mathematically equivalent, 
but (\ref{def_TNN_normed}) has better numerical stability during the training process.
From Figure \ref{TNNstructure} and numerical tests, we can find the parameters 
for each rank of TNN are
correlated by the FNN, which guarantee the stability of the TNN-based machine 
learning methods in some sense. 
This is also an important difference from the tensor finite element methods.

In order to show the reasonableness of TNN, we now introduce the approximation 
property from \cite{WangJinXie}, which is based on 
the isomorphism relation between $H^m(\Omega_1\times\cdots\times\Omega_d)$
and the tensor product space $H^m(\Omega_1)\otimes\cdots\otimes H^m(\Omega_d)$. 
The process of approximating the function $f(x)\in H^m(\Omega_1\times\cdots\times\Omega_d)$
by the TNN defined as (\ref{def_TNN_normed}) can be regarded as searching 
for a correlated CP decomposition structure
to approximate $f(x)$ in the space $H^m(\Omega_1)\otimes\cdots\otimes H^m(\Omega_d)$
with the rank being not greater than $p$. 
The following approximation result to the functions in the space 
$H^m(\Omega_1\times\cdots\times\Omega_d)$ under the sense of $H^m$-norm 
is proved in \cite{WangJinXie}.
\begin{theorem}\cite[Theorem 1]{WangJinXie}\label{theorem_approximation}
Assume that each $\Omega_i$ is an interval in $\mathbb R$ for $i=1, \cdots, d$, $\Omega=\Omega_1\times\cdots\times\Omega_d$,
and the function $f(x)\in H^m(\Omega)$. Then for any tolerance $\varepsilon>0$, there exist a
positive integer $p$ and the corresponding TNN defined by (\ref{def_TNN_normed})
such that the following approximation property holds
\begin{equation}\label{eq:L2_app}
\|f(x)-\Psi(x;\theta)\|_{H^m(\Omega)}<\varepsilon.
\end{equation}
\end{theorem}

Although there is no general result to give the relationship between the hyperparameter $p$ 
and error bounds, we also provided an estimate of the rank $p$ 
under a smoothness assumption. For easy understanding, we focus on the periodic setting with 
$I^d=I\times I\times\cdots\times I=[0,2\pi]^d$ and the approximations 
property of TNN to the functions in the linear space which is defined with Fourier basis.
Note that similar approximation results of TNN can be extended to the non-periodic functions. 
For each variable $x_i\in[0,2\pi]$, let us define the one-dimensional 
Fourier basis $\{\varphi_{k_i}(x_i):= \frac{1}{\sqrt{2\pi}}e^{{\rm i}k_ix_i},k_i\in\mathbb Z\}$ 
and classify functions via the decay of their Fourier coefficients. 
Further denote multi-index $k=(k_1,\cdots,k_d)\in\mathbb Z^d$ and $x=(x_1,\cdots,x_d)\in I^d$. 
Then the $d$-dimensional Fourier basis can be built with the tensor product way 
\begin{eqnarray}
\varphi_k(x)\coloneqq\prod_{i=1}^d\varphi_{k_i}(x_i)=(2\pi)^{-d/2}e^{{\rm i}k\cdot x}.
\end{eqnarray}
We denote
\begin{eqnarray}
\lambda_{\rm mix}(k)\coloneqq\prod_{i=1}^d\left(1+|k_i|\right)\ \ \ 
{\rm and}\ \ \ \lambda_{\rm iso}(k)\coloneqq 1+\sum_{i=1}^d|k_i|.
\end{eqnarray}
Now, for $-\infty<t,\ell<\infty$, we define the space $H_{\rm mix}^{t,\ell}(I^d)$ 
as follows (cf. \cite{GriebelHamaekers,Knapek}) 
{\footnotesize 
\begin{eqnarray}
H_{\rm mix}^{t,\ell}(I^d)=\left\{u(x)=\sum_{k\in\mathbb Z^d}c_k\varphi_k(x):
\|u\|_{H_{\rm mix}^{t,\ell}(I^d)}
=\left(\sum_{k\in\mathbb Z^d}\lambda_{\rm mix}(k)^{2t}\cdot\lambda_{\rm iso}(k)^{2\ell}\cdot|c_k|^2\right)^{1/2}<\infty \right\}. 
\end{eqnarray}}
Note that the parameter $\ell$ governs the isotropic smoothness, 
whereas $t$ governs the mixed smoothness. The space $H_{\rm mix}^{t,\ell}(I^d)$ 
gives a quite flexible framework for the study of problems in Sobolev spaces.
See \cite{GriebelHamaekers,GriebelKnapek,Knapek} for more information 
on the space $H_{\rm mix}^{t,\ell}(I^d)$. 
Thus, \cite{WangJinXie} gives the following comprehensive error estimate for TNN. 
\begin{theorem}\label{theorem_aprrox_rate}
Assume function $f(x)\in H_{\rm mix}^{t,\ell}(I^d)$, $t>0$ and $m>\ell$. 
Then there exists a TNN $\Psi(x;\theta)$ 
defined by (\ref{def_TNN_normed}) such that the following approximation property holds
\begin{eqnarray}
\|f(x)-\Psi(x;\theta)\|_{H^m(I^d)}\leq C(d)\cdot p^{-(\ell-m+t)}
\cdot\|u\|_{H_{\rm mix}^{t,\ell}(I^d)},
\end{eqnarray}
where $C(d)\leq c\cdot d^2\cdot0.97515^d$,  $c$ is independent of $d$ 
and each subnetwork of TNN is a FNN which is built by using $\sin(x)$ 
as the activation function and one hidden layer with $2p$ neurons. 
\end{theorem}

From the approximation result in Theorem \ref{theorem_aprrox_rate}, 
when the target function belongs to $H_{\rm mix}^{t,\ell}(\Omega)$,  
there exists a TNN with $p\sim\mathcal O(\varepsilon^{-(m-\ell-t)})$ 
such that the accuracy is $\varepsilon$.
For more details about the approximation properties of TNN, please refer to \cite{WangJinXie}. 

\section{Tensor neural network based machine learning method 
}\label{homo-TNN}
In this section, the TNN based machine learning method will be introduced 
for solving the elliptic multiscale problem. We first introduce the solving procedure 
and then show the way to compute the high dimensional integration included in the loss functions.

\subsection{The solving procedure of TNN-based method}
In order to design the machine learning method for MPDEs, we execute the homogenization 
method described in Section \ref{Section_MPDEs}.  Then the machine learning method 
can be described as follows: 
\begin{enumerate}
\item Solve the cell problem (\ref{cell-prob}) using TNN based machine learning method.  

\item Compute the homogenized coefficients based on the solutions of the cell problem. 
Iterate the process until we obtain the approximation for the homogenized coefficient $A_0^*$. 

\item 
Solve the homogenized equation (\ref{homo-homogenization}) by using TNN based machine learning method.
\end{enumerate}
For easy understanding, we take the case of $K=1$ as the example. 
In the first step, in order to solve the cell problems (\ref{cell-prob}), 
we build the TNN function $\Psi_j(\mathbf x, \mathbf y, \Theta)$ 
described in Section \ref{Section_TNN} as the approximation to the functions 
in $L^2\left(D; H_{\mathrm{per}}^1(Y)\right)$.
In order to make $\Psi_j(\mathbf x, \mathbf y, \Theta)$ 
satisfy the periodic boundary condition on $\mathbf{y}$, 
we choose the sine function as the activation function, and fix all entries in the 
weight matrix of the first layer for each FNN on the corresponding dimension 
of $\mathbf{y}$ as the multiples of $2\pi$. 
Meanwhile, the corresponding biases are registered as trainable parameters, 
and so are the remaining parameters for $\mathbf{y}$ and all the parameters 
for $\mathbf{x}$. Under these settings, $\Psi_j(\mathbf x, \mathbf y, \Theta)$ 
satisfies the periodic boundary condition on $\mathbf{y}$ automatically.
In addition, since the integral of the solution of equation (\ref{cell-prob}) is $0$ 
with respect to $\mathbf{y}$, we construct a new tensor neural network function 
based on $\Psi_j(\mathbf x, \mathbf y, \Theta)$ as
\begin{eqnarray}\label{hat-psi}
\widehat{\Psi}_j(\mathbf x, \mathbf y, \Theta) = \Psi_j(\mathbf x, \mathbf y, \Theta) 
- \int_Y \Psi_j(\mathbf x, \mathbf y, \Theta)d\mathbf y ,\ \ \forall j = 1,\cdots,d,
\end{eqnarray}
which satisfies the following condition
\begin{eqnarray}\label{Integration_Constraint}
\int_Y \widehat{\Psi}_j(\mathbf x, \mathbf y, \Theta)d\mathbf y = 0 ,
\ \ \ \forall \mathbf x\in \Omega,\ \ \forall j = 1,\cdots,d.
\end{eqnarray}
The TNN function $\widehat{\Psi}_j(\mathbf x, \mathbf y, \Theta)$ now satisfies 
the constraints of the solution for the equation (\ref{cell-prob}), 
and therefore can be used to approximate the exact solution $\chi_j(\mathbf x, \mathbf y)$. Notice that it is usually 
tricky to make the trial functions satisfy the integration constraint (\ref{Integration_Constraint}) of the cell problem (\ref{cell-prob}), 
since the accuracy of the high dimensional integration is always hard to be guaranteed.  
However, due to the special structure of TNN, the constraint here is easy and natural to deal with for TNN-based method.
Next, we define the corresponding loss function as follows
\begin{eqnarray}\label{loss_1step}
L_j(\theta) = \left\|\frac{\partial}{\partial y_i}\left(a_{ik}(\mathbf x, \mathbf y)
\frac{\partial \widehat{\Psi}_j(\mathbf x, \mathbf y,\theta)}{\partial y_k}\right) 
+ \frac{\partial a_{ij}(\mathbf x, \mathbf y)}{\partial y_i}\right\|_{L^2(\Omega;Y)}, 
\ \ j = 1, \cdots, d.
\end{eqnarray}
Here, we use Einstein notation. Then by solving the following optimization problems
\begin{eqnarray}\label{approx_opt_1step}
\min_{\theta\in \Theta} \ \ L_j(\theta) , \ \ j = 1, \cdots, d,
\end{eqnarray}
we obtain the optimal parameters $\theta^*_j$, and the corresponding 
$\widehat{\Psi}_j(\mathbf x, \mathbf y, \theta^*_j)$, $1\leq j\leq d$, 
are the numerical approximations to the cell problems (\ref{cell-prob}). 

In the second step, we compute the homogenized coefficient as follows
\begin{eqnarray}\label{coeff}
a_{i j}^{\text{TNN}}(\mathbf{x})
=\int_Y\left(a_{i j}(\mathbf{x},\mathbf{y})+a_{i k}(\mathbf{x},\mathbf{y})
\frac{\widehat{\Psi}_j(\mathbf{x},\mathbf{y},\theta^*_j)}{\partial y_k}\right)
d\mathbf{y}.
\end{eqnarray}
Based on the tensor product structure of $\widehat{\Psi}_j(\mathbf{x},\mathbf{y},\theta^*_j)$, 
we are able to compute the homogenized coefficient defined in (\ref{coeff}) with high accuracy.

In the last step, we build another TNN function $\Phi(\mathbf x, \Theta)$ 
as the approximation to the solution $u_0(\mathbf x)\in H_0^1\left(D\right)$ 
of the homogenized equation (\ref{homo-homogenization}).  
In order to satisfy the Dirichlet boundary condition of the problem (\ref{homo-homogenization}), 
the TNN function $\widehat{\Phi}(\mathbf x, \Theta)$ is defined as 
\begin{eqnarray}\label{hat-phi}
\widehat{\Phi}(\mathbf x, \Theta) = g(\mathbf x){\Phi}(\mathbf x, \Theta),
\end{eqnarray}
where $g(\mathbf x) = 0, \forall \mathbf x\in  \partial \Omega$ and $ g(\mathbf x) \neq 0$, 
$\forall \mathbf x \in \Omega \backslash \partial \Omega$.
Then we define the corresponding loss function as follows
\begin{eqnarray}\label{loss_3step}
L(\theta) = \left\|\frac{\partial}{\partial x_i}\left(a_{ik}^{\text{TNN}}(\mathbf x)
\frac{\partial \widehat{\Phi}(\mathbf x,\theta)}{\partial x_k}\right)
+f(\mathbf x)\right\|_{L^2(\Omega)},
\end{eqnarray}
and solve the following optimization problem
\begin{eqnarray}\label{approx_opt_3step}
\min_{\theta\in \Theta} \ \ L(\theta)
\end{eqnarray}
to obtain the optimal parameters $\theta^*$. The corresponding 
$\widehat{\Phi}(\mathbf x, \theta^*)$ is the numerical solution 
for the homogenized problem (\ref{homo-homogenization}). 
Through the tricks mentioned above, the TNN functions (\ref{hat-psi}) 
and (\ref{hat-phi}) satisfy the corresponding boundary condition 
automatically. Therefore, (\ref{approx_opt_1step}) and (\ref{approx_opt_3step}), 
what we need to solve, are all unconstrained optimization problems. 

Different from the normal FNN-based machine learning method,  
where the Monte Carlo integration is usually the indispensable option, 
the quadrature scheme with fixed quadrature points can be used 
to do the numerical integrations in this paper. 
Fortunately, based on TNN structure in the loss function (\ref{approx_opt_1step}) 
and (\ref{approx_opt_3step}), Theorem 3 in \cite{WangJinXie}
shows that these numerical integrations here does not encounter ``curse of dimensionality''
since the computational work can be bounded by the polynomial scale of dimension $d$. 
Due to the high accuracy of the high dimensional integrations with Gauss points, 
all the computational process can be done with high accuracy, 
which is the core reason why the TNN-based method achieves high accuracy.

\subsection{Quadrature scheme for TNN}\label{Section_Integration}
In this subsection, we introduce the quadrature scheme for computing the integrations 
in the loss function (\ref{approx_opt_1step}), (\ref{approx_opt_3step}) and the homogenized 
coefficient (\ref{coeff}), which involves TNN functions and tensor-product-type multiscale coefficients.
 As for the general case, please refer to \cite{WangJinXie}, where the method to compute the numerical 
 integrations of polynomial composite functions of TNN and their derivatives are designed. 
 In the loss function (\ref{approx_opt_1step})
\begin{equation*}
\begin{aligned}
&\ \ \ \ \ L_j(\theta) = \left\|\frac{\partial}{\partial y_i}\left(a_{ik}(\mathbf x, \mathbf y)
\frac{\partial \widehat{\Psi}_j(\mathbf x, \mathbf y,\theta)}{\partial y_k}\right) 
+ \frac{\partial a_{ij}(\mathbf x, \mathbf y)}{\partial y_i}\right\|_{L^2(\Omega;Y)} \\
&= \left(\int_{\Omega \times Y} \left(\sum_{i=1}^d\frac{\partial}{\partial y_i}\left(\sum_{k=1}^da_{ik}(\mathbf x, \mathbf y)
\frac{\partial \widehat{\Psi}_j(\mathbf x, \mathbf y,\theta)}{\partial y_k}\right) 
+ \sum_{i=1}^d\frac{\partial a_{ij}(\mathbf x, \mathbf y)}{\partial y_i}\right)^2
d\mathbf x d \mathbf y \right)^{\frac{1}{2}} \\
&= \left(\int_{\Omega \times Y} \left(\sum_{i=1}^d \sum_{k=1}^d 
\left(\frac{\partial a_{ik}(\mathbf x, \mathbf y)}{\partial y_i} 
\frac{\partial \widehat{\Psi}_j(\mathbf x, \mathbf y,\theta)}{\partial y_k} + a_{ik}(\mathbf x, \mathbf y)\frac{\partial^2 \widehat{\Psi}_j(\mathbf x, \mathbf y,\theta)}{\partial y_k \partial y_i} \right) 
+ \sum_{i=1}^d\frac{\partial a_{ij}(\mathbf x, \mathbf y)}{\partial y_i}\right)^2 
d\mathbf x d \mathbf y \right)^{\frac{1}{2}}, 
\end{aligned}
\end{equation*}
we need to compute a $(K+1)d$-dimensional integration. 
With the tensor-product-type structure of TNN and the multiscale coefficients, 
we are able to compute the integration with polynomial scale computational 
complexity and high accuracy by transforming the high dimensional integration 
into the multiplication of $(K+1)d$ one-dimensional integrations. 
As an example, we consider the case of $K=1, d=2$, and compute 
one of the terms of the $L^2$ square expansion, i.e., 
\begin{equation}\label{integration_square_expansion}
\int_{\Omega  \times Y} \left(\frac{\partial a_{11}(\mathbf x, \mathbf y)}{\partial y_1} 
\frac{\partial \widehat{\Psi}_j(\mathbf x, \mathbf y,\theta)}{\partial y_1}\right)^2 d\mathbf x d \mathbf y.
\end{equation}
We make two assumptions here. One is that $\Omega$ can be decomposed into the tensor product of two one-dimensional 
set, i.e., $\Omega = \Omega_1 \times \Omega_2$, where $\Omega_1$ and $\Omega_2$ can be a bounded interval $(a_i,b_i)$, 
the whole line $(-\infty,+\infty)$ or the half line $(a_i,+\infty)$. 
Furthermore, we assume that the entries of the multiscale coefficients $A$ have tensor-product structure, 
i.e., they are in the form of
\begin{equation*}
a_{ik}(\mathbf{x}, \mathbf{y}) = \sum_{e=1}^{p_{ik}} \alpha_{ik,e} \phi_{ik,1,e}(x_1)\phi_{ik,2,e}(x_2)\phi_{ik,3,e}(y_1)\phi_{ik,4,e}(y_2).
\end{equation*}
Note that there are some other tricks to deal with the non-tensor-product-type case, 
the basic idea in which is to use some tensor-product-type basis functions, 
such as Fourier basis, polynomial basis, and even TNN functions in the form of (\ref{def_TNN_normed}), 
to interpolate the coefficient. We will cover these tricks in an upcoming article.
Recall that by (\ref{hat-psi}) and (\ref{def_TNN_normed}), $\widehat{\Psi}_j$ is also tensor-product-type. 
For the sake of better understanding, we omit some indices here and set
\begin{equation*}
\begin{aligned}
a_{11} := \sum_{e=1}^{p} \alpha_{1,e} \prod_{i=1}^4 \phi_{1,i,e}(x_i),\ \ \ \ \ 
\widehat{\Psi}_j := \sum_{\ell=1}^q \alpha_{2,\ell} \prod_{i=1}^4 \phi_{2,i,\ell}(x_i).
\end{aligned}
\end{equation*}
Here, we set $x_3 := y_1$, $x_4 := y_2$, $\Omega_3 := [0,1]$ and $\Omega_4 := [0,1]$ for the sake of notational simplicity.
Hence, (\ref{integration_square_expansion}) can be computed by
\begin{eqnarray}\label{intergration_square_tensor_form}
&& \int_{\Omega  \times Y} \left(\frac{\partial a_{11}(\mathbf x, \mathbf y)}{\partial y_1} 
\frac{\partial \widehat{\Psi}_j(\mathbf x, \mathbf y,\theta)}{\partial y_1}\right)^2 
d\mathbf x d \mathbf y \nonumber\\
&&= \int_{\Omega_1 \times \Omega_2 \times \Omega_3 \times \Omega_4} 
\left(\sum_{e=1}^{p} \alpha_{1,e} 
\prod_{i \neq 3} \phi_{1,i,e}(x_i) \frac{\partial \phi_{1,3,e}(x_3)}{\partial x_3} \cdot 
\sum_{\ell=1}^q \alpha_{2,\ell} \prod_{i \neq 3} \phi_{2,i,\ell}(x_i)
\frac{\partial \phi_{2,3,\ell}(x_3)}{\partial x_3} \right)^2 \nonumber\\
&& \quad \qquad \qquad \qquad \qquad \qquad \times dx_1 dx_2 dx_3 dx_4\nonumber\\
&&=\sum_{e=1}^p \sum_{j=1}^p \sum_{\ell=1}^q \sum_{k=1}^q \alpha_{1,e} \alpha_{1,j} \alpha_{2,\ell} \alpha_{2,k} \cdot \prod_{i \neq 3} 
\int_{\Omega_i} \phi_{1,i,e}(x_i)\phi_{1,i,j}(x_i)\phi_{2,i,\ell}(x_i)\phi_{2,i,k}(x_i) dx_i\nonumber \\
&& \quad \qquad \qquad \qquad \qquad \qquad \cdot \int_{\Omega_3} 
\frac{\partial \phi_{1,3,e}(x_3)}{\partial x_3} \frac{\partial \phi_{1,3,j}(x_3)}{\partial x_3} 
\frac{\partial \phi_{2,3,\ell}(x_3)}{\partial x_3} \frac{\partial \phi_{2,3,k}(x_3)}{\partial x_3} dx_3.
\end{eqnarray}
The expansion (\ref{intergration_square_tensor_form}) gives the hint to design the efficient numerical
scheme to compute the high dimensional integration (\ref{integration_square_expansion}).
Without loss of generality, for the $i$-th one dimensional domain $\Omega_i$,
we choose $N_i$ Gauss points $\{x_i^{(n_i)}\}_{n_i=1}^{N_i}$
and the corresponding weights $\{w_i^{(n_i)}\}_{n_i=1}^{N_i}$ to compute the one
dimensional integrations included in  (\ref{intergration_square_tensor_form}), where $i=1, \cdots, d$.   
Then the splitting numerical integration scheme for (\ref{integration_square_expansion}) can be given 
as follows:
\begin{equation*}
\begin{aligned}
& \int_{\Omega  \times Y} \left(\frac{\partial a_{11}(\mathbf x, \mathbf y)}{\partial y_1} \frac{\partial \widehat{\Psi}_j(\mathbf x, \mathbf y,\theta)}{\partial y_1}\right)^2 d\mathbf x d \mathbf y \\
\approx & \sum_{e=1}^p \sum_{j=1}^p \sum_{\ell=1}^q \sum_{k=1}^q \alpha_{1,e} \alpha_{1,j} \alpha_{2,\ell} \alpha_{2,k} \cdot \prod_{i \neq 3} \sum_{n_i=1}^{N_i}w_i^{n_i} \phi_{1,i,e}(x_i^{n_i})\phi_{1,i,j}(x_i^{n_i})\phi_{2,i,\ell}(x_i^{n_i})\phi_{2,i,k}(x_i^{n_i}) \\
& \quad \qquad \qquad \qquad \qquad \qquad \cdot \sum_{n_3=1}^{N_3} w_3^{n_3} \frac{\partial \phi_{1,3,e}(x_3^{n_3})}{\partial x_3} \frac{\partial \phi_{1,3,j}(x_3^{n_3})}{\partial x_3} \frac{\partial \phi_{2,3,\ell}(x_3^{n_3})}{\partial x_3} \frac{\partial \phi_{2,3,k}(x_3^{n_3})}{\partial x_3}, 
\end{aligned}
\end{equation*}
and therefore, the computational work of the high-dimensional integration of (\ref{integration_square_expansion}) is reduced to the polynomial scale of dimension $d$. 
Following the same procedure, we are able to compute all the terms of the $L^2$ 
square expansion of loss function (\ref{approx_opt_1step}), (\ref{approx_opt_3step}) 
and the homogenized coefficient (\ref{coeff}) with high accuracy and acceptable computational complexity. 
Solving the minimization problems using optimization algorithms to find the optimized parameters 
of the neural networks becomes much more reasonable and numerically stable 
when the accuracy of computing loss functions is guaranteed, 
which is the initial motivation for designing the TNN structure.

\section{Numerical results}\label{Section_Numerical}
In this section, we provide several examples to validate the efficiency and accuracy 
of the proposed TNN-based machine learning method for MPDEs. 
For comparison, we solve the multiscale problem directly with the finite element method 
on the fine enough mesh, and the corresponding finite element solution 
$u_{\varepsilon}^{\text{FEM}}$ acts as the reference to check the accuracy.

\subsection{One-dimensional elliptic two-scale problems}
In this subsection, we solve a one-dimensional two-scale problem (i.e. $d = 1$, $K = 1$) 
with the TNN-based machine learning method.   
Let $\Omega$ be the interval $[0,\pi]$ and $Y$ be the unit interval $[0,1]$. 
We consider the following problem
\begin{eqnarray}\label{ex_1D}
\left\{
\begin{array}{rcl}
-\frac{\partial }{\partial x} \left(A(x,\frac{x}{\varepsilon})
\frac{\partial u_{\varepsilon}}{\partial x}\right) = f, & x\in \Omega,& \\
u_{\varepsilon}(0) = u_{\varepsilon}(\pi) = 0,&&
\end{array}
\right.
\end{eqnarray}
where the multiscale coefficient $A(x,y) = 0.5\sin(2\pi y) + \sin(x) + 2$, 
and the function on the right hand side $ f(x) = \sin(x)$.  

In order to solve  (\ref{cell-prob}) and (\ref{homo-homogenization}), 
we need to calculate the integrations in the loss functions (\ref{loss_1step}) 
and (\ref{loss_3step}). For this aim,  the one dimensional interval 
$\Omega=[0, \pi]$ and $Y=[0,1]$ are decomposed 
into $40$ equal subintervals and $16$ Gauss points are used on each subinterval.

\begin{figure}[htb!]
\centering
\includegraphics[width=1\textwidth]{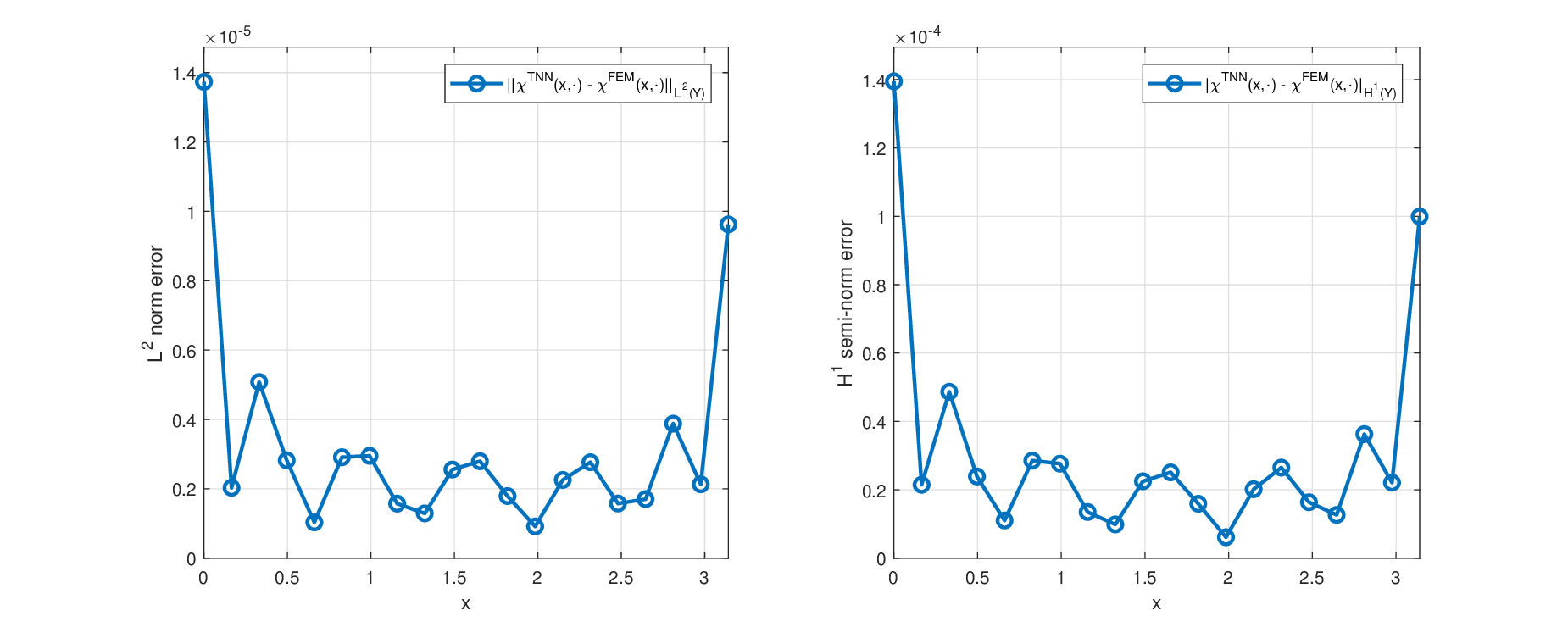}
\caption{The $L^2$ error and $H^1$ semi-norm error of 
$\chi^{\text{TNN}}(x_i, \cdot)$ at the selected points $\{x_i\}, 
i=1,\cdots, 20$. On the selected points $x_i$, the maximum $L^2$ 
norm error is 1.3736e-05, the maximum $H^1$ semi-norm error is 1.3949e-04, 
the maximum relative $L^2$ norm error is 4.7948e-04, 
the maximum relative $H^1$ semi-norm error is 7.7027e-04.}\label{fig:comp-ex1-1D}
\end{figure}

For the problem (\ref{cell-prob}), we choose the rank $p=10$ and the 
Adam optimizer \cite{KingmaAdam} with learning rate $0.01$ is adopted 
in the first $30000$ steps and then the LBFGS with learning rate $0.1$ 
is used for the subsequent $20000$ steps. Each subnetwork of TNN is chosen 
as the FNN with two hidden layers and each hidden layer has $20$ neurons, 
see Figure \ref{TNNstructure}. In order to validate the accuracy of 
TNN based machine learning method for $\chi_1({x}, {y})$ of the problem (\ref{cell-prob}), 
we select $20$ uniformly distributed spatial points  $\{x_i\}, i=1,\cdots,20$ on $\Omega$ 
and solve the second-order elliptic equation (\ref{cell-prob}) 
with respect to $y$ for each $x_i$ with the finite element method.  
The corresponding finite element approximations $\chi^{\text{FEM}}(x_i, \cdot)$ 
act as the reference solution to measure the accuracy of the TNN 
solution $\chi^{\text{TNN}}(x_i, \cdot)$ at the selected points. 
Figure \ref{fig:comp-ex1-1D} shows the corresponding error estimates 
of the TNN solution $\chi^{\text{TNN}}(x_i, \cdot)$ 
in the $L^2$ norm and $H^1$ semi-norm, respectively,  at the selected points.
It can be easily found that the TNN based machine learning method 
can reach high accuracy for the cell problem (\ref{cell-prob}).

\begin{figure}[htb!]
\centering
\includegraphics[width=0.4\textwidth]{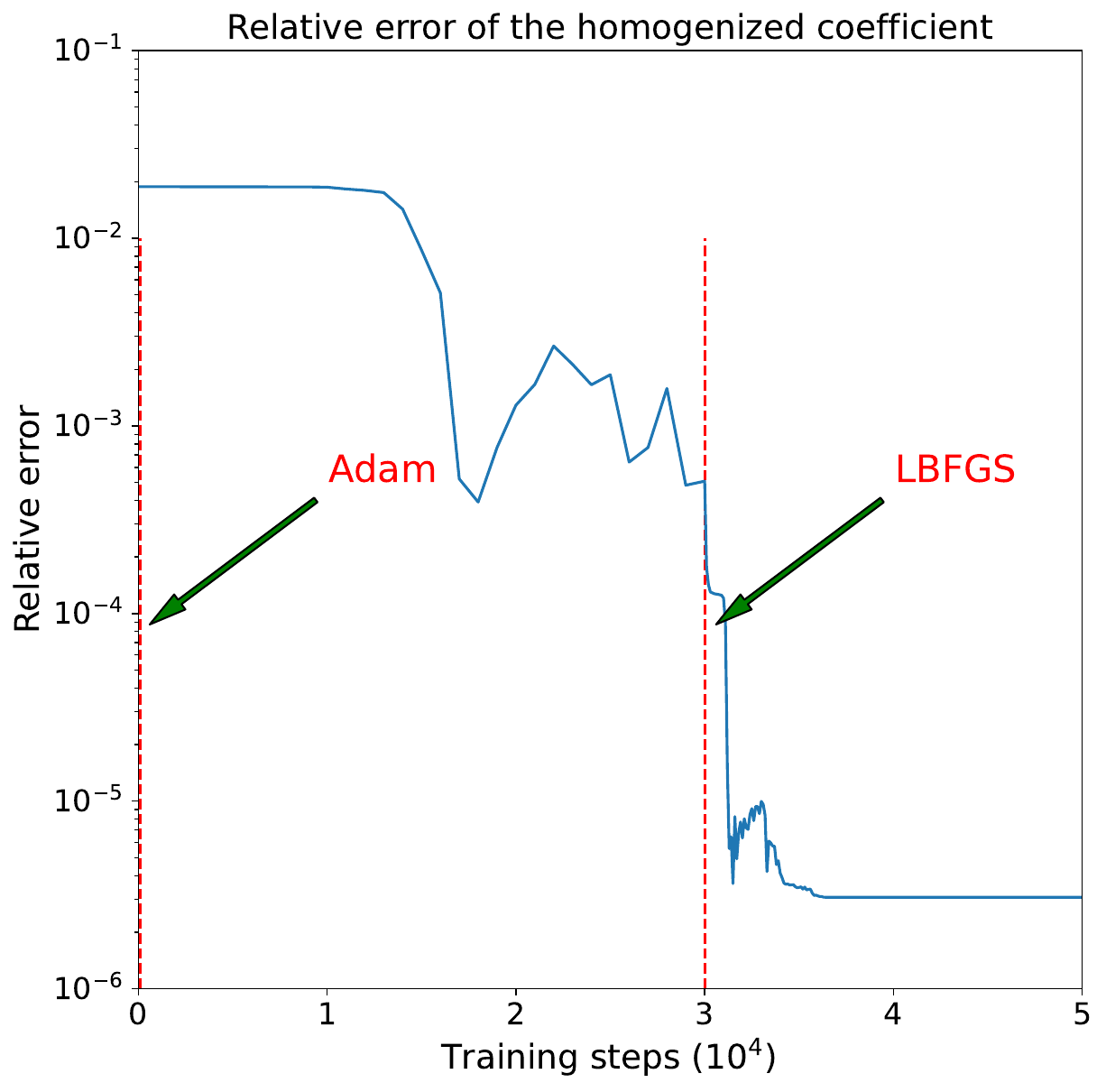}
\caption{Relative errors of the homogenized coefficient during the 
training process for problem (\ref{cell-prob}).}
\label{fig:1D_homo_coeff_training}
\end{figure}

\begin{figure}[htb!]
\centering
\includegraphics[width=1\textwidth]{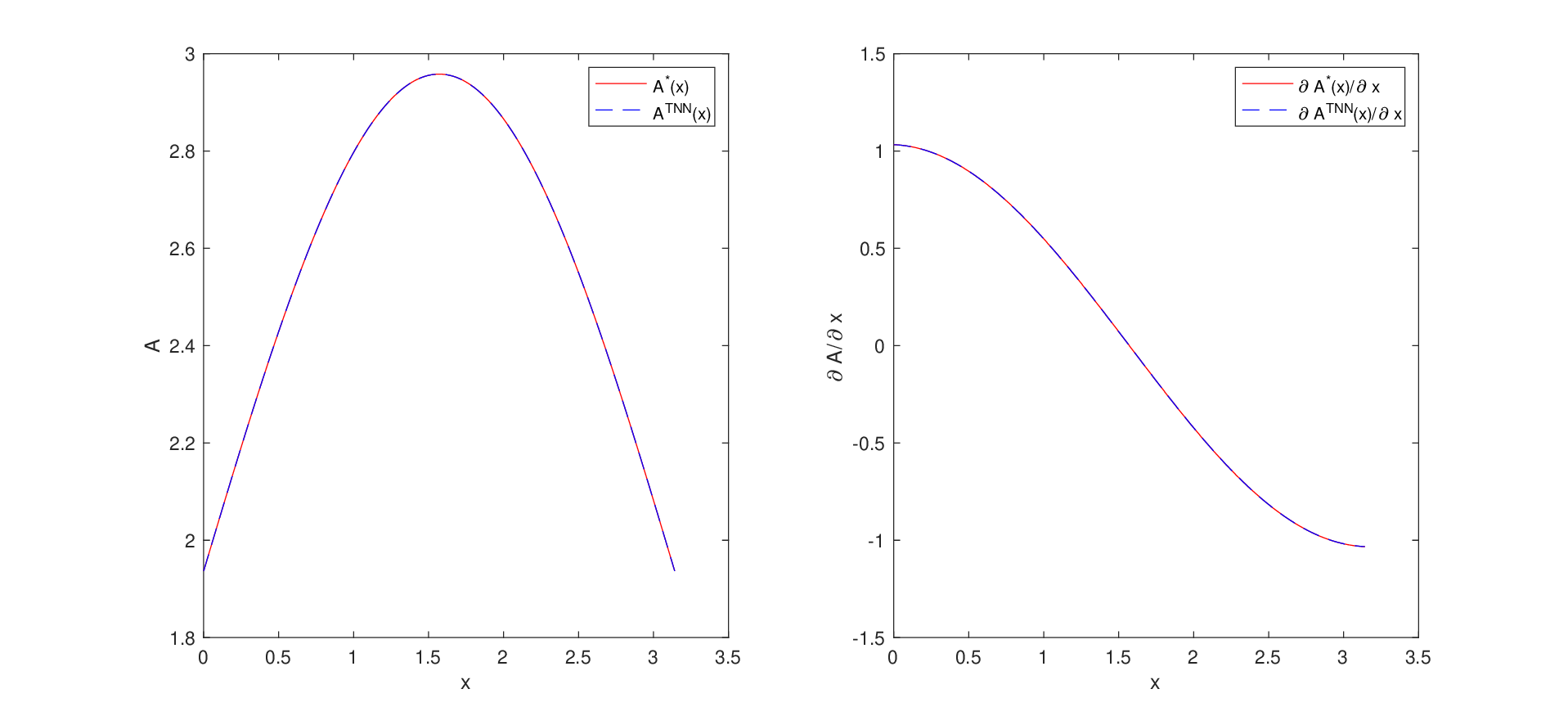}
\caption{The figures of $A^{\text{TNN}}(x)$, 
$\frac{\partial A^{\text{TNN}}(x)}{\partial x}$ 
and their exact ones for the problem (\ref{ex_1D}), 
where the exact homogenized coefficient $A^*(x)$ is defined by 
formula (\ref{homo-coefficient-analytical}) 
and $A^{\text{TNN}}(x)$ is the approximate coefficient obtained by the TNN method. 
Left: the homogenized coefficients $A^*(x)$ and  $A^{\text{TNN}}(x)$. 
Right: the homogenized coefficient $\frac{\partial A^*(x)}{\partial x}$ 
and $\frac{\partial A^{\text{TNN}}(x)}{\partial x}$.}\label{fig:1D_Astar}
\end{figure}

Then, the homogenized coefficient $A^{\text{TNN}}(x)$ can be calculated 
from $\chi^{\text{TNN}}(x, y)$ by formula (\ref{coeff}).  
The analytical homogenized coefficient $A^*(x)$ defined in (\ref{homo-coefficient-analytical}) 
is used for comparison. The final relative $L^2$ norm error estimate of 
the homogenized coefficient calculated by the TNN method is $\|A^{\text{TNN}}-A^*\|_{L^2(\Omega)}/\|A^*\|_{L^2(\Omega)} = 3.0678\text{e-06}$. 
Figure \ref{fig:1D_homo_coeff_training} shows the relative error 
of $A^{\text{TNN}}$ during the training process.
The final approximations $A^{\text{TNN}}(x)$ and $\frac{\partial A^{\text{TNN}}(x)}{\partial x}$, 
the analytical homogenized coefficient are demonstrated in Figure \ref{fig:1D_Astar} for comparison.

For problem (\ref{homo-homogenization}), we first build a TNN with the same structure 
as that in problem (\ref{cell-prob}) and then use the Adam optimizer with learning rate of 0.01
for the first 20000 training steps and LBFGS with learning rate of 0.1 for the subsequent 30000 
steps. In order to investigate the accuracy of TNN solution $u_0^{\text{TNN}}$,  
the finite element solution of the equation (\ref{homo-homogenization}) 
with the accurate homogenized coefficient $A^*(x)$ in (\ref{homo-coefficient-analytical}) 
is used as the reference. The absolute and relative error estimates in the $L^2$ norm 
and $H^1$ semi-norm for $u_0^{\text{TNN}}$ are as below:
\begin{eqnarray*}
\|u_0^{\text{TNN}}-u_0^{\text{FEM}}\|_{L^2(\Omega)} = 3.3264\text{e{-05}},&& \frac{\|u_0^{\text{TNN}}-u_0^{\text{FEM}}\|_{L^2(\Omega)}}{\|u_0^{\text{FEM}}\|_{L^2(\Omega)}} = 6.2106\text{e{-05}}, \nonumber\\
|u_0^{\text{TNN}}-u_0^{\text{FEM}}|_{H^1(\Omega)} = 1.6429\text{e{-04}},&& \frac{|u_0^{\text{TNN}}-u_0^{\text{FEM}}|_{H^1(\Omega)}}{|u_0^{\text{FEM}}|_{H^1(\Omega)}} 
= 3.0497\text{e{-04}}.
\end{eqnarray*}

\begin{figure}[htb!]
\centering
\includegraphics[width=1\textwidth]{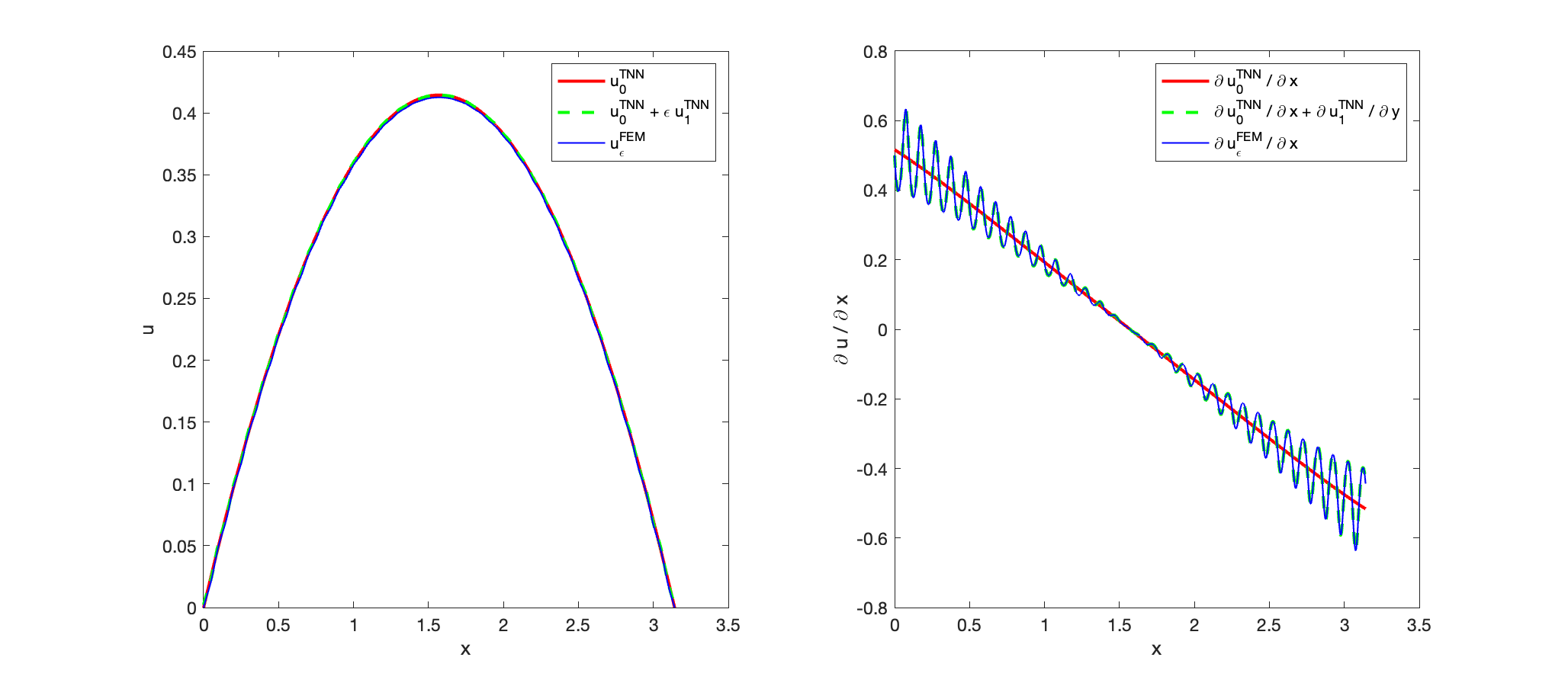}
\caption{Left: the figures of $u_0^{\text{TNN}}$, 
$u_0^{\text{TNN}} + \varepsilon u_1^{\text{TNN}} $ and $  
u_{\varepsilon}^{\text{FEM}}$ of the problem (\ref{ex_1D}) 
for the case of $\varepsilon = 1/10$. 
Right: the figures of $\frac{\partial u_0^{\text{TNN}}}{\partial x}$, 
$\frac{\partial u_0^{TNN}}{\partial x} + \frac{\partial u_1^{TNN}}{\partial y}$ 
and $\frac{\partial u_{\varepsilon}^{\text{FEM}}}{\partial x}$ 
of the problem (\ref{ex_1D})  for the case of $\varepsilon = 1/10$. 
Here, the final $H^1$ semi-norm absolute error 
$|u_{0}^{\text{TNN}} + \varepsilon u_1^{\text{TNN}} 
- u_{\varepsilon}^{\text{FEM}}|_{H^1(\Omega)} = 2.699\text{e{-04}}$.}
\label{fig:1D_ue}
\end{figure}

\begin{figure}[htb]
\centering
\includegraphics[width=0.5\textwidth]{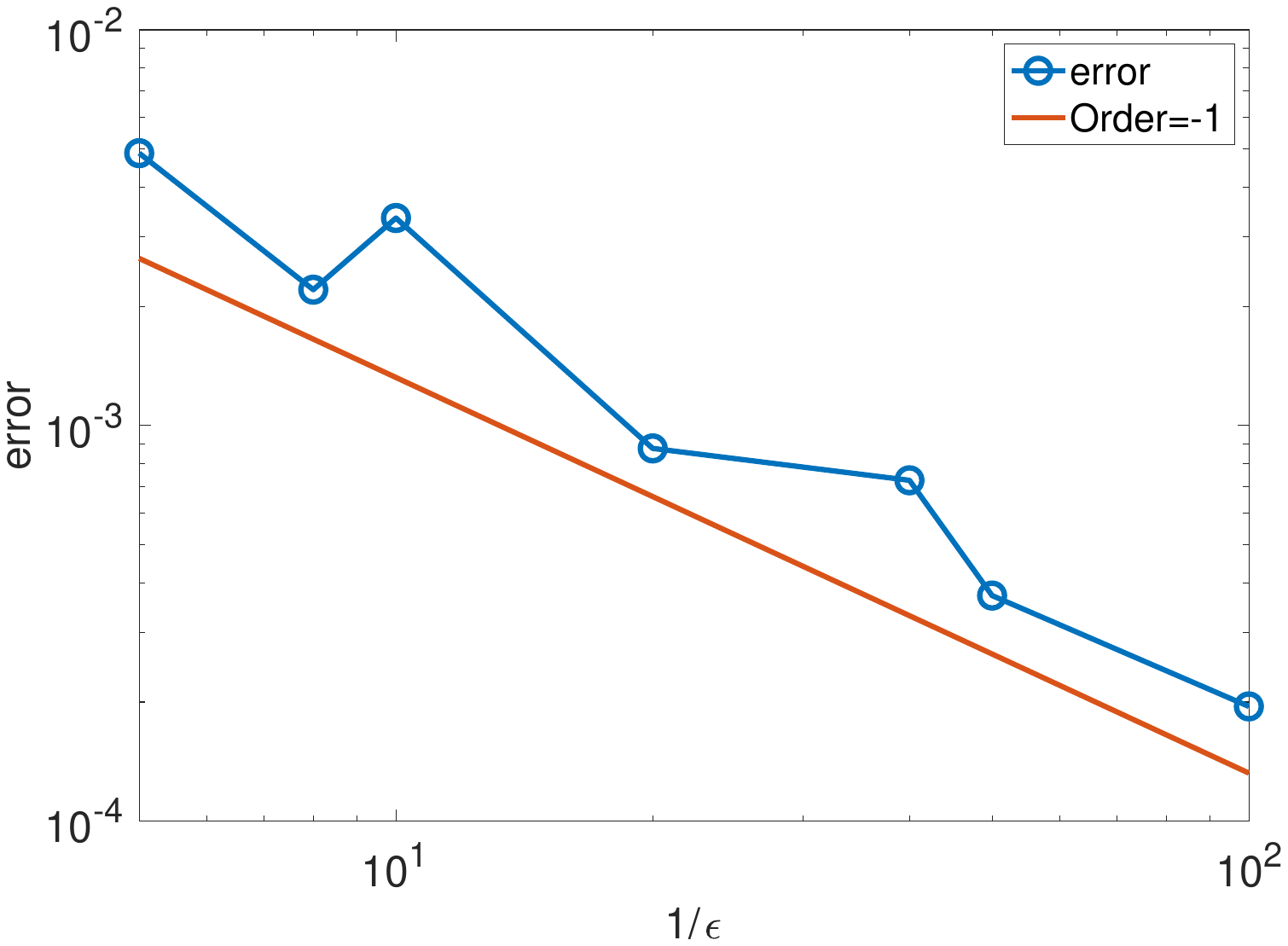}
\caption{The $L^2$ norm errors of TNN solutions $u_{\varepsilon}^{\text{TNN}} 
= u_{0}^{\text{TNN}} + \varepsilon u_1^{\text{TNN}}$ for different $\varepsilon$ of the 
problem (\ref{ex_1D})}\label{fig:1D_homoerror}
\end{figure}

Finally, in order to check the convergence results (\ref{Convergence_L2}) 
and (\ref{Convergence_H1}), we investigate the accuracy of the multiscale solution 
$u_0^{\text{TNN}} + \varepsilon u_1^{\text{TNN}}$ by TNN based machine learning method 
under different  $\varepsilon$.  
For different $\varepsilon$, the multiscale problem is solved directly by the 
finite element method on a sufficiently fine mesh and the corresponding solution 
$u_{\varepsilon}^{\text{FEM}}$ is adopted as the reference. 
The TNN solutions and finite element solution for the case of 
$\varepsilon = 1/10$ are shown in Figure \ref{fig:1D_ue}. 
Furthermore, Table \ref{tab:1D_error_eps} and Figure \ref{fig:1D_homoerror} show
the error estimates $\|u_{0}^{\text{TNN}} + \varepsilon u_1^{\text{TNN}} - u_{\varepsilon}^{\text{FEM}}\|_{L^2(\Omega)}$ 
in the $L^2$ norm of problem (\ref{ex_1D}) 
for different $\varepsilon$. The results in Figure \ref{fig:1D_homoerror} validate the  
convergence results of the multiscale expansion (\ref{Convergence_L2}), 
and the TNN based method can reach high accuracy with absolute error $1.9492\text{e{-04}}$ 
and relative error $3.6408\text{e{-04}}$ in the $L^2$ norm for the case of $\varepsilon = 0.01$. 
The TNN based method is expected to have higher accuracy for smaller $\varepsilon$.
\begin{table}[htb]
\centering
\caption{The $L^2$ norm absolute errors 
$\|u_{0}^{\text{TNN}} + \varepsilon u_1^{\text{TNN}} - u_{\varepsilon}^{\text{FEM}}\|_{L^2(\Omega)}$ 
and relative errors $\|u_{0}^{\text{TNN}} 
+ \varepsilon u_1^{\text{TNN}} - u_{\varepsilon}^{\text{FEM}}\|_{L^2(\Omega)}/ \|u_{\varepsilon}^{\text{FEM}}\|_{L^2(\Omega)}$ 
between the TNN solution and the reference solution 
of problem (\ref{ex_1D}) for different $\varepsilon$.}
\label{tab:1D_error_eps}
\begin{tabular}{||c|c|c|c|c|c|c||}
\hline
$\varepsilon$   & $\frac{1}{5}$& $\frac{1}{8}$ & $\frac{1}{10}$ & $\frac{1}{20}$ & $\frac{1}{50}$ & $\frac{1}{100}$\\
\hline
Absolute errors & 4.8789e-03 & 2.2040e-03 &3.3453e-03 & 8.7515e-04 & 3.7163e-04 &1.9497e-04\\
\hline
Relative errors & 9.1813e-03 & 4.1212e-03 &6.2812e-03 & 1.6349e-03 & 6.9388e-04 &3.6408e-04\\
\hline 
\end{tabular}
\end{table}

It is also worth mentioning that, the problem (\ref{homo-homogenization}) 
can also be solved by using TNN-based a posteriori error adaptive algorithm described 
in \cite{TNN-posterior}. Using the same TNN structure and the same training parameters 
setting, we do obtain a more accurate $\widehat{u}_0^{\text{TNN}}$. 
The absolute and relative errors in $L^2$ norm, $H^1$ semi-norm of 
$\widehat{u}_0^{\text{TNN}}$ are as below:
\begin{eqnarray*}
\|\widehat{u}_0^{\text{TNN}}-u_0^{\text{FEM}}\|_{L^2(\Omega)} = 3.3370\text{e{-07}},&&
\frac{\|\widehat{u}_0^{\text{TNN}}-u_0^{\text{FEM}}\|_{L^2(\Omega)}}{\|u_0^{\text{FEM}}\|_{L^2(\Omega)}}
= 6.2304\text{e{-07}},\\
|\widehat{u}_0^{\text{TNN}}-u_0^{\text{FEM}}|_{H^1(\Omega)} = 2.6831\text{e{-06}}, && \frac{|\widehat{u}_0^{\text{TNN}}-u_0^{\text{FEM}}|_{H^1(\Omega)}}{|u_0^{\text{FEM}}|_{H^1(\Omega)}} 
= 4.9807\text{e{-06}},
\end{eqnarray*}
which show that $\widehat{u}_0^{\text{TNN}}$ has better accuracy than that of $u_0^{\text{TNN}}$. 
Since the multiscale expansion error is independent of the accuracy of the numerical approximations to 
the limit problems, we do not continue checking the improved accuracy $\|\widehat{u}_{0}^{\text{TNN}} 
+ \varepsilon \widehat{u}_1^{\text{TNN}} - u_{\varepsilon}^{\text{FEM}}\|_{L^2(\Omega)}$. 
But the TNN-based a posteriori error adaptive algorithm in \cite{TNN-posterior} can be used to treat 
smaller $\varepsilon$ with much smaller accuracy 
$\|\widehat{u}_{0}^{\text{TNN}} + \varepsilon \widehat{u}_1^{\text{TNN}} 
- u_{\varepsilon}^{\text{FEM}}\|_{L^2(\Omega)}$ to balance smaller multiscale expansion accuracy.

\subsection{Two-dimensional elliptic two-scale problems}
This subsection is devoted to investigating the TNN-based machine learning method 
for solving two-dimensional two-scale problems (i.e. $d = 2$, $K = 1$). 
We use the same quadrature scheme as that of the example in the last subsection, 
i.e., all intervals in all 4 dimensions are decomposed into 40 equal subintervals 
and 16 Gauss points are chosen on each subinterval. 

First, we consider the following second order elliptic equation from \cite{LeungLinZhang}.
Let $\Omega$ be the unit square $[0,1]^2$ and corresponding $Y$ be the unit square $[0,1]^2$. 
Here, we are concerned with the following MPDE 
\begin{eqnarray}\label{ex_2D_1}
\left\{
\begin{array}{rcl}
-\frac{\partial }{\partial x_i} \left(A(\frac{\mathbf{x}}{\varepsilon})
\frac{\partial u_{\varepsilon}}{\partial x_i}\right) = f, & \mathbf{x}\in \Omega,& \\
u_{\varepsilon}(\mathbf{x}) = 0,&\mathbf{x}\in\partial\Omega,&
\end{array}
\right.
\end{eqnarray}
where the multiscale coefficient $A(\mathbf{y}) = 2 + \sin(2\pi y_1)\cos(2\pi y_2)$,  
and the function $ f(\mathbf{x}) = \sin(x_1) + \cos(x_2).$ 

\begin{figure}[htb!]
\centering
\subfigure{
\includegraphics[width=0.48\textwidth]{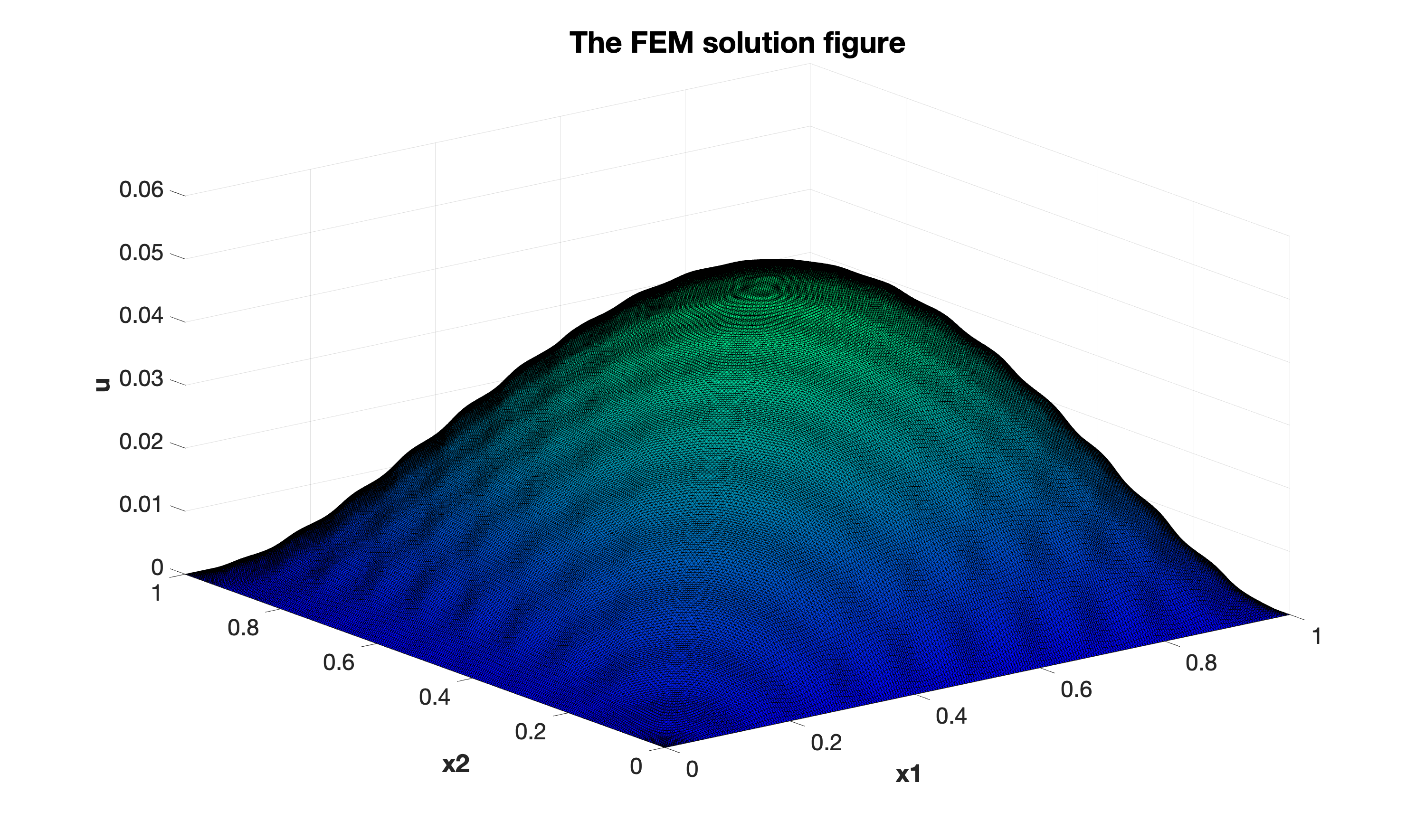}
}
\subfigure{
\includegraphics[width=0.48\textwidth]{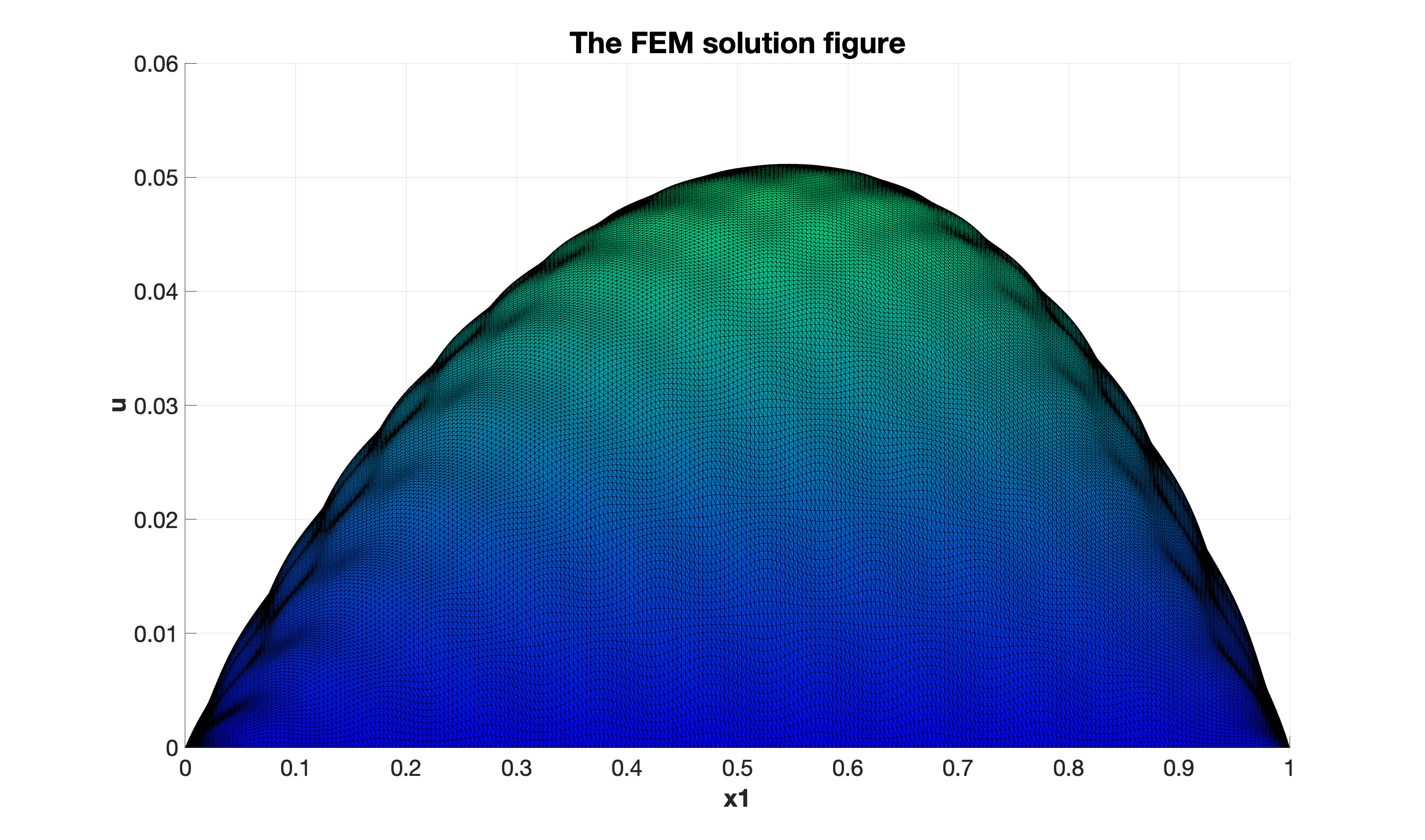}
}

\subfigure{
\includegraphics[width=0.48\textwidth]{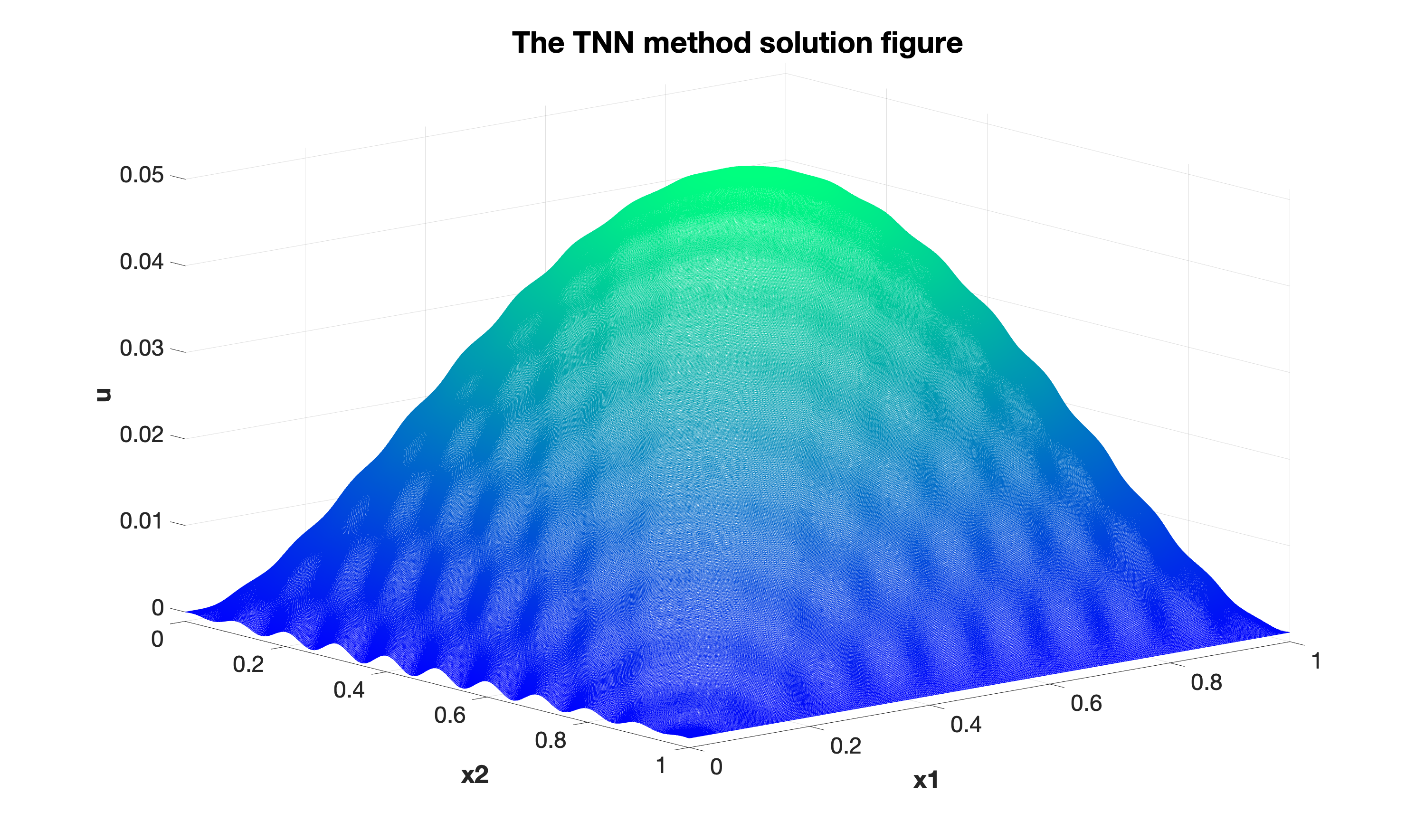}
}
\subfigure{
\includegraphics[width=0.48\textwidth]{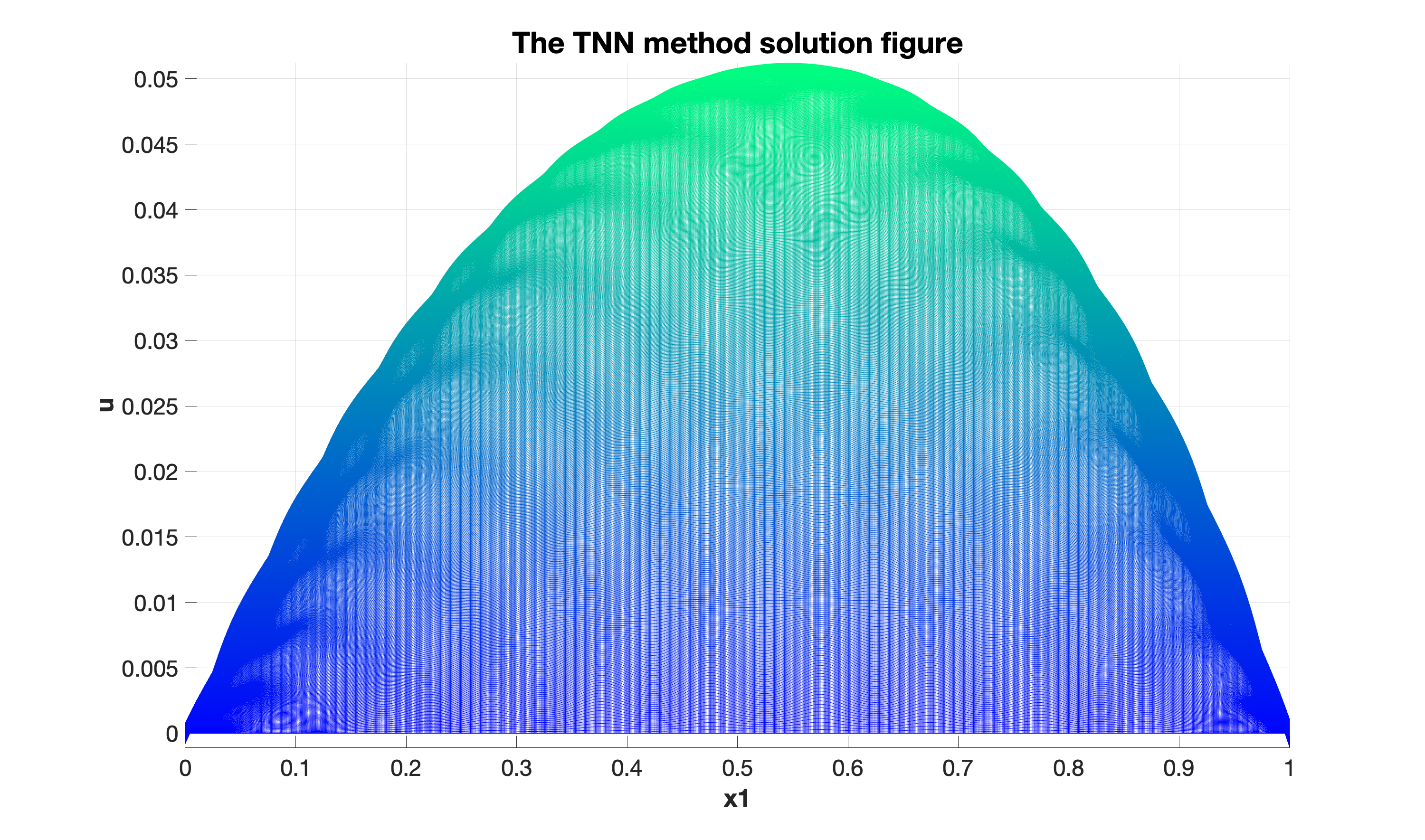}
}
\caption{The solution of finite element method and TNN method for 
problem (\ref{ex_2D_1}) with $\varepsilon = 1/10$. 
Up-Left: $u_{\varepsilon}^{\text{FEM}}$. Up-Right: X-Z view 
of $u_{\varepsilon}^{\text{FEM}}$. Down-Left: $u_{0}^{\text{TNN}} 
+ \varepsilon u_{1}^{\text{TNN}} $. Down-Right: X-Z view of 
$u_{0}^{\text{TNN}} + \varepsilon u_{1}^{\text{TNN}} $.} 
\label{fig:Ay_2scale_e10} 
\end{figure}

The homogenization method described in Section \ref{homo-TNN} is adopted to solve this problem. 
For problem (\ref{cell-prob}), we build the TNN with $p = 20$ and 
each subnetwork being a FNN with two hidden layers and $20$ neurons 
for each hidden layer. In order to train the TNN, the Adam optimizer with learning 
rate 0.01 is adopted for the first 20000 steps and then the LBFGS with learning 
rate 0.1 for the subsequent 30000 steps. 
The next step is to obtain the homogenized coefficient matrix $A^*$ by the formula (\ref{coeff}) 
based on the TNN solutions of problem (\ref{cell-prob}).
For problem (\ref{homo-homogenization}), we use the Adam optimizer with learning rate 0.01 
for the first 20000 steps and then the LBFGS with learning rate 0.1 for the subsequent 30000 steps 
to train the TNN which has the same structure as that in problem (\ref{cell-prob}). 
Similarly, the reference solution is obtained by using the finite element method 
on a sufficiently fine mesh for different $\varepsilon$. 
The finite element solution and the TNN solution is shown in Figure \ref{fig:Ay_2scale_e10} 
for the case of $\varepsilon = 1/10$.
\begin{table}[htb!]
\centering
\caption{The $L^2$ norm absolute errors $\|u_{0}^{\text{TNN}} + \varepsilon u_1^{\text{TNN}} 
- u_{\varepsilon}^{\text{FEM}}\|_{L^2(\Omega)}$ and relative errors $\|u_{0}^{\text{TNN}} 
+ \varepsilon u_1^{\text{TNN}} - u_{\varepsilon}^{\text{FEM}}\|_{L^2(\Omega)}/ 
\|u_{\varepsilon}^{\text{FEM}}\|_{L^2(\Omega)}$ 
between the TNN solution and the reference solution of the problem (\ref{ex_2D_1}) 
for different $\varepsilon$.}\label{tab:2D_Ay_error_eps}
\begin{tabular}{||c|c|c|c|c|c||}
\hline \hline
$\varepsilon$ & $\frac{1}{5}$ & $\frac{1}{8}$ & $\frac{1}{10}$ & $\frac{1}{30}$ &  $\frac{1}{100}$\\
\hline
Absolute errors& 2.2304e-04 & 1.1601e-04 & 8.6779e-05 & 2.7003e-05 & 9.9181e-06\\
\hline
Relative errors& 7.8946e-03 & 4.1032e-03 & 3.0684e-03 & 9.5398e-04 & 3.5027e-04\\
\hline \hline
\end{tabular}
\end{table}
\begin{figure}[htb!]
\centering
\includegraphics[width=0.5\textwidth]{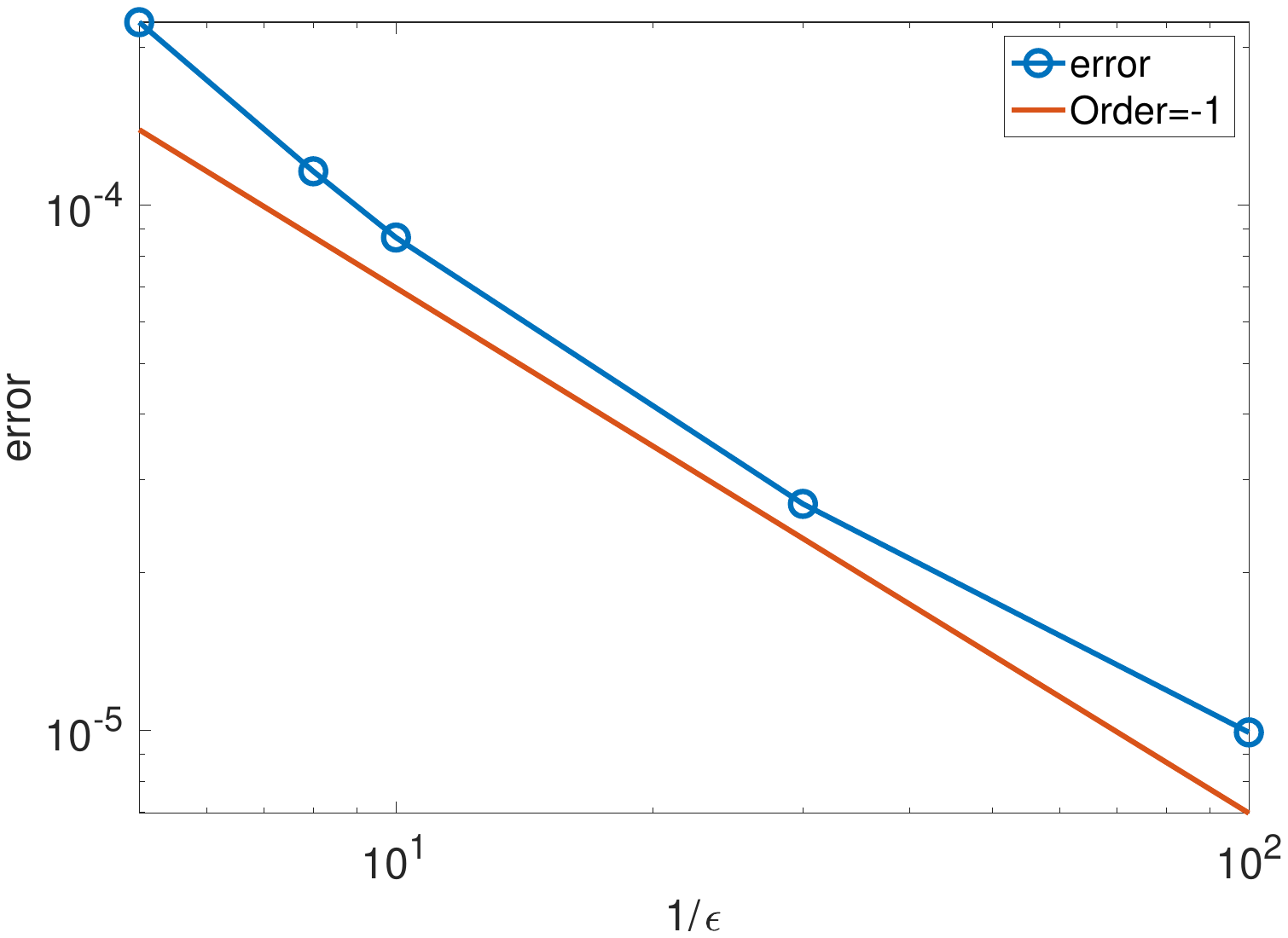}
\caption{The $L^2$ norm errors of TNN solutions $u_{\varepsilon} 
= u_0^{\text{TNN}} + \varepsilon u_1^{\text{TNN}}$ for 
different $\varepsilon$ of the problem (\ref{ex_2D_1})}\label{fig:2D_homoerror_Ay}
\end{figure}

Table \ref{tab:2D_Ay_error_eps} and Figure \ref{fig:2D_homoerror_Ay} show 
the error results $\|u_{0}^{\text{TNN}} + \varepsilon u_1^{\text{TNN}} - u_{\varepsilon}^{\text{FEM}}\|_{L^2(\Omega)}$ 
for the problem (\ref{ex_2D_1}). 
We can also find that the numerical results in Figure \ref{fig:2D_homoerror_Ay} 
validate the theoretical convergence results (\ref{Convergence_L2}) 
for the multiscale expansion method. 
The TNN based method reaches the high accuracy  $9.9181\text{e{-06}}$ 
for the absolute error and $3.5027\text{e{-04}}$ for the relative error in the 
$L^2$ norm for the case of $\varepsilon = 0.01$.

In the second example, the coefficient $A^{\varepsilon}$ has both the fast and slow variables.
Let $\Omega$ be $[0,\pi]^2$ and corresponding $Y$ be the unit square $[0,1]^2$. 
And we consider the following problem
\begin{eqnarray}\label{ex_2D_2}
\left\{
\begin{array}{rcl}
-\frac{\partial }{\partial x_i} \left(A(\mathbf{x},\frac{\mathbf{x}}{\varepsilon})\frac{\partial u_{\varepsilon}}{\partial x_i}\right) = f, & \mathbf{x}\in \Omega,& \\
u_{\varepsilon}(\mathbf{x}) = 0,&\mathbf{x}\in\partial\Omega,&
\end{array}
\right.
\end{eqnarray}
where the multiscale coefficient $A(\mathbf{x},\mathbf{y}) = 0.5 \sin(2 \pi y_1) \sin(2 \pi y_2) 
+ \sin(x_1) + \sin(x_2) + 3$ and the function $f(\mathbf{x}) = \sin(x_1) \sin(x_2)$.  

\begin{figure}[htb!]
\centering
\includegraphics[width=0.5\textwidth]{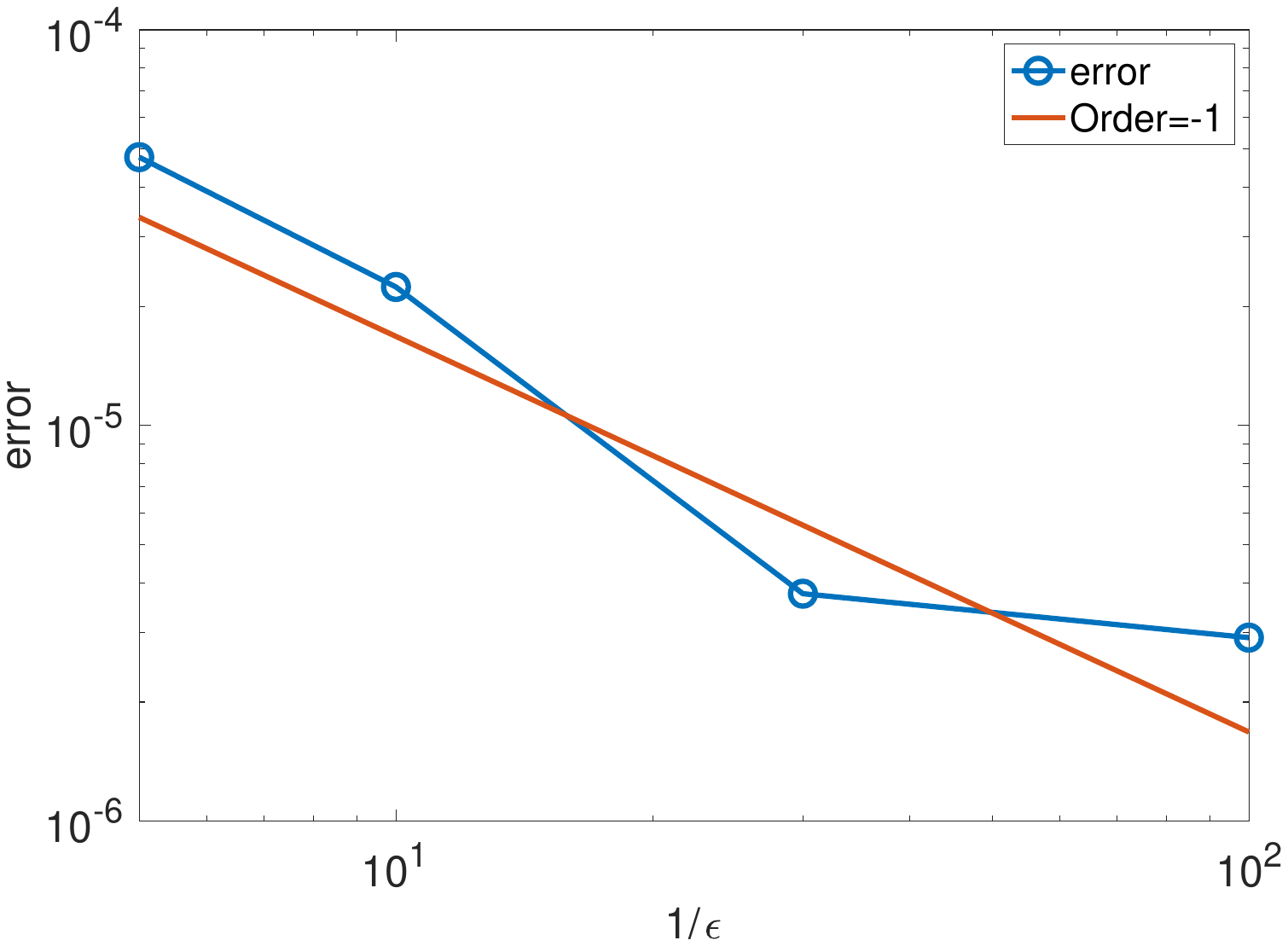}
\caption{The $L^2$ norm errors of TNN solutions $u_{\varepsilon} = u_0^{\text{TNN}} 
+ \varepsilon u_1^{\text{TNN}}$ of the problem (\ref{ex_2D_2}) for different $\varepsilon$.}
\label{fig:2D_homoerror_Axy}
\end{figure}

With the three-steps procedure, the same TNN structure and parameters 
setting as that for the problem (\ref{ex_2D_1}) are adopted here.
We also investigate the accuracy of multiscale solutions by TNN based method 
for different $\varepsilon$ and the reference solution is also obtained by using 
the finite element method on the fine enough mesh.
Table \ref{tab:2D_Axy_error_eps} and Figure \ref{fig:2D_homoerror_Axy} 
show the $L^2$ error $\|u_{0}^{\text{TNN}} + \varepsilon u_1^{\text{TNN}} 
- u_{\varepsilon}^{\text{FEM}}\|_{L^2}$ for the problem (\ref{ex_2D_2}). 
\begin{table}[htb!]
\caption{The $L^2$ norm absolute errors $\|u_{0}^{\text{TNN}} + \varepsilon u_1^{\text{TNN}} - u_{\varepsilon}^{\text{FEM}}\|_{L^2(\Omega)}$ 
and relative errors $\|u_{0}^{\text{TNN}} + \varepsilon u_1^{\text{TNN}} 
- u_{\varepsilon}^{\text{FEM}}\|_{L^2(\Omega)}/ \|u_{\varepsilon}^{\text{FEM}}\|_{L^2(\Omega)}$ 
between the TNN solution and the reference solution of the problem (\ref{ex_2D_2}) 
for different $\varepsilon$.}\label{tab:2D_Axy_error_eps}
\centering
\begin{tabular}{||c|c|c|c|c||}
\hline \hline
$\varepsilon$ &$\frac{1}{5}$ & $\frac{1}{10}$ & $\frac{1}{30}$ & $\frac{1}{100}$\\
\hline
Absolute errors & 4.7672e-05 & 2.2396e-05 & 3.7570e-06 & 2.9091e-06\\
\hline
Relative errors & 2.5750e-04 & 1.2097e-04 & 2.0292e-05 & 1.5713e-05\\
\hline \hline
\end{tabular}
\end{table}
Similarly, the numerical results in Figure \ref{fig:2D_homoerror_Axy} 
validate the theoretical convergence results of the homogenization 
expansion (\ref{Convergence_L2}). The TNN based method obtains 
high accuracy $2.9091\text{e{-06}}$ for the absolute error 
and  $1.5713\text{e{-05}}$ for the relative error in the $L^2$ 
norm for the case of $\varepsilon = 0.01$.

When the coefficient $A$ has both fast and slow variables, the deduced 
cell equation (\ref{cell-prob}) is a high dimensional problem defined 
in the space $\mathbb R^{2d}$, which is very difficult to solve 
by the traditional numerical methods. In this case, 
the TNN based machine learning method can also solve 
this high dimensional cell problem with high accuracy, 
which is the main contribution of this paper. 

\subsection{One-dimensional elliptic three-scale problems}
In this subsection, we consider a one-dimensional elliptic equation 
with three scales (i.e. $d = 1$, $K = 2$). Let $\Omega$ be 
the interval $[0, \pi]$ and $Y$ be the unit interval $[0,1]$. 
We consider the following problem
\begin{eqnarray}\label{ex_1D_3scale}
\left\{
\begin{array}{rcl}
-\frac{\partial }{\partial x} \left(A(x,\frac{x}{\varepsilon},
\frac{x}{\varepsilon^2})\frac{\partial u_{\varepsilon}}{\partial x}\right) 
= f, & x\in \Omega,& \\
u_{\varepsilon}(0) = u_{\varepsilon}(\pi) = 0,&&
\end{array}
\right.
\end{eqnarray}
where the multiscale coefficient $A(x,y_1,y_2) = 0.5\sin(2\pi y_2) 
+ \sin(2 \pi y_1) + \sin(x) + 3$, and the function on the right 
hand side $ f(x) = \sin(x)$. 

We now describe the detailed procedure to solve the MPDE (\ref{ex_1D_3scale}). 
For the quadrature scheme, $\Omega$ and $Y$ are both decomposed into 100 equal subintervals 
and 8 Gauss points are chosen on each subinterval. 
For the problem (\ref{cell-prob}) with coefficient $A$ 
and divergence with respect to $y_2$, we choose 
the rank $p=20$ and the Adam optimizer with learning rate $0.01$ is 
adopted for the first $20000$ steps and then the LBFGS optimizer with 
learning rate 0.1 for the subsequent $20000$ steps. Each subnetwork 
of TNN is chosen as the FNN with two hidden layers and each hidden 
layer has 20 neurons. We then compute the approximate coefficient 
$A_1^{\text{TNN}}(x, y_1)$ and use the analytical coefficient $A_1^*(x, y_1)$ 
obtained by formula (\ref{homogenized-matrix-analytical}) 
for comparison. Following up, the problem (\ref{cell-prob}) 
with coefficient $A_1^{\text{TNN}}$ and divergence with respect to $y_1$ 
is solved by the TNN based machine learning method. 
For this aim, we choose the same TNN structure and parameters setting as the first cell 
problem expect with a learning rate of 0.3 for the LBFGS optimizer.

\begin{figure}[htb!]
\centering
\subfigure[Relative errors of $A_1^{\text{TNN}}$]
{\includegraphics[width=0.4\textwidth]{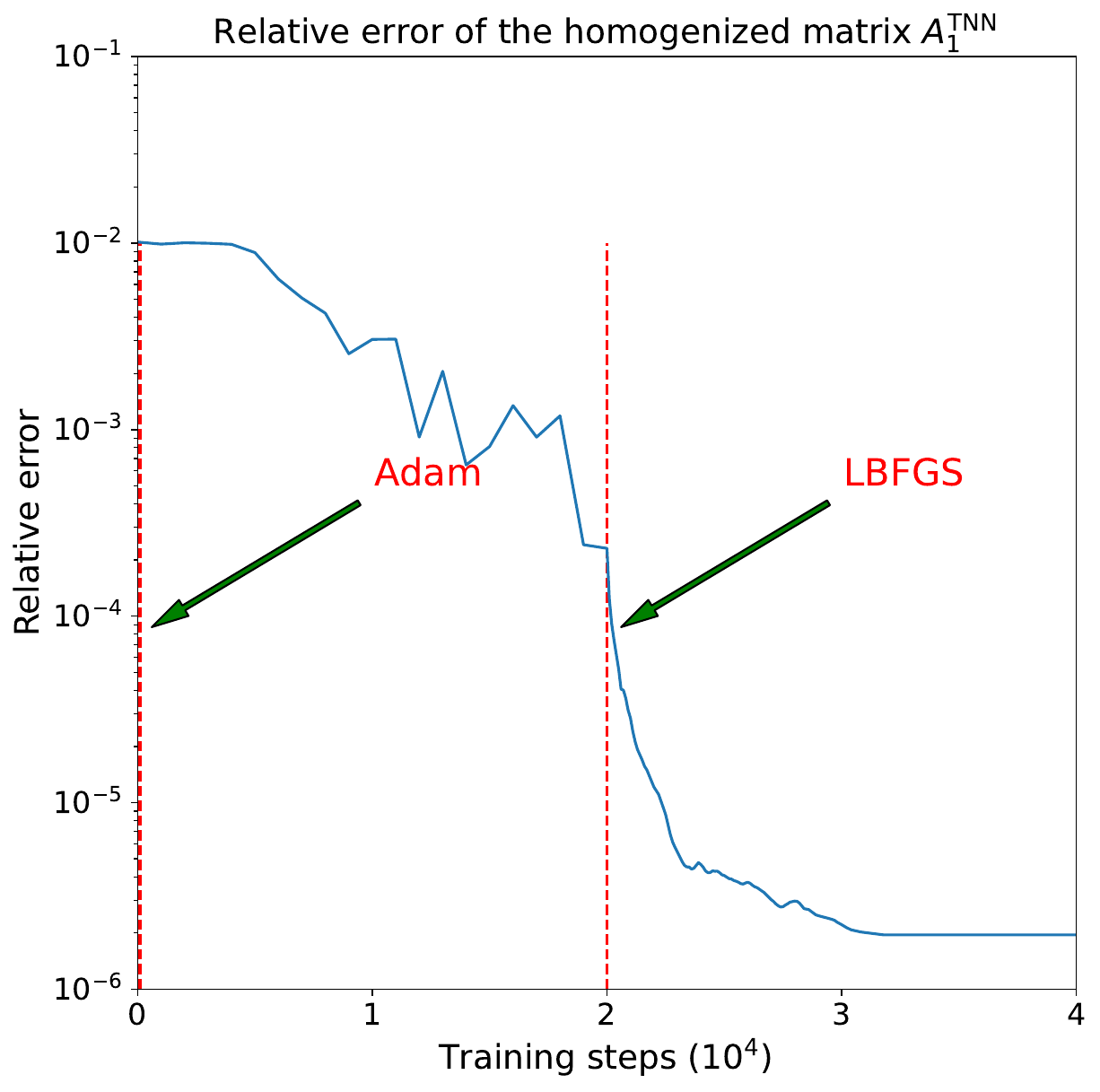}}
\subfigure[Relative errors of $A_0^{\text{TNN}}$]
{\includegraphics[width=0.4\textwidth]{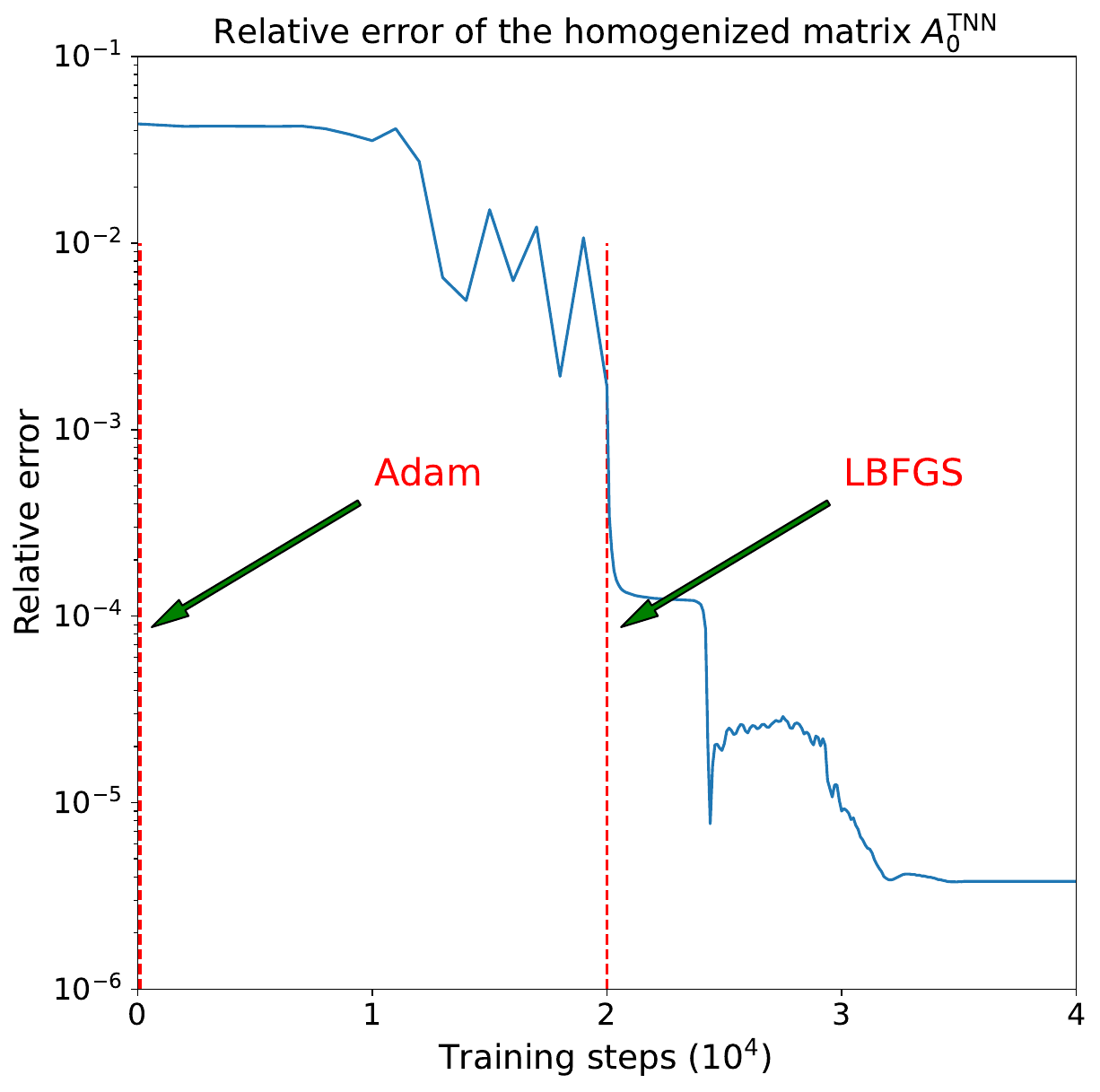}}
\caption{Relative errors of the homogenized matrix $A_1^{\text{TNN}}$ 
and $A_0^{\text{TNN}}$ during the training process.}
\label{fig:1D_3scales_homo_coeff_training}
\end{figure}

Again we compute the approximate coefficient $A_0^{\text{TNN}}(x)$ 
and use the analytical coefficient $A_0^*(x)$ obtained by 
formula (\ref{homo-coefficient-analytical}) for comparison. 
Figure \ref{fig:1D_3scales_homo_coeff_training} 
shows the relative errors of $A_1^{\text{TNN}}$ and $A_0^{\text{TNN}}$ 
during the training process. The final relative $ L^2$ norm errors 
of these two homogenized matrices are 
$\|A_1^{\text{TNN}}-A_1^*\|_{L^2(\Omega \times Y_1)}/\|A_1^*\|_{L^2(\Omega \times Y_1)} 
= 1.9551\text{e-06}$ and $\|A_0^{\text{TNN}}-A_0^*\|_{L^2(\Omega)}/\|A_0^*\|_{L^2(\Omega)} 
= 3.7828\text{e-06}$.

For the final homogenized equation (\ref{homo-homogenization}) with homogenized 
coefficient $A_0^{\text{TNN}}$, we adopt the same TNN structure and parameters 
setting as that of the first cell problem. To validate the effectiveness of 
TNN for this step, the equation (\ref{homo-homogenization}) with the analytical 
homogenized coefficient $A_0^*(x)$ is solved by the finite element method 
on fine enough mesh and the corresponding solution $u_0^{\text{FEM}}$ 
acts as the reference. 
The absolute and relative $L^2$ norm error of the TNN solution $u_0^{\text{TNN}}$ 
are as follows:
\begin{eqnarray*}
\|u_0^{\text{TNN}}-u_0^{\text{FEM}}\|_{L^2(\Omega)} = 2.7074\text{e{-07}},&& \frac{\|u_0^{\text{TNN}}-u_0^{\text{FEM}}\|_{L^2(\Omega)}}{\|u_0^{\text{FEM}}\|_{L^2(\Omega)}} 
= 6.9219\text{e{-07}}. \nonumber
\end{eqnarray*}

Finally, we conduct a back substitute process using (\ref{uk-uk-1}) and investigate 
the accuracy of the multiscale solution $u_0^{\text{TNN}} + \varepsilon u_1^{\text{TNN}} 
+ \varepsilon^2 u_2^{\text{TNN}}$ for $\varepsilon = 1/5$, $1/10$ and $1/30$, respectively. 
Similarly, we use the finite element method on a sufficiently fine mesh to 
solve (\ref{ex_1D_3scale}) directly to obtain the reference solution 
$u_\varepsilon^{\text{FEM}}$. The absolute and relative $L^2$ norm 
errors of $u_0^{\text{TNN}} + \varepsilon u_1^{\text{TNN}} 
+ \varepsilon^2 u_2^{\text{TNN}}$ are shown in Table \ref{tab:1D_3scale_error_eps}. 
As expected, the accuracy is higher for smaller $\varepsilon$, 
and it is safe to conclude that TNN based machine learning method is 
able to obtain sufficiently accurate solution for smaller $\varepsilon$.
\begin{table}[htb!]
\caption{The $L^2$ norm absolute errors $\|u_{0}^{\text{TNN}} 
+ \varepsilon u_1^{\text{TNN}} + \varepsilon^2 u_2^{\text{TNN}} 
- u_{\varepsilon}^{\text{FEM}}\|_{L^2(\Omega)}$ 
and relative errors $\|u_{0}^{\text{TNN}} + \varepsilon u_1^{\text{TNN}} 
+ \varepsilon^2 u_2^{\text{TNN}} 
- u_{\varepsilon}^{\text{FEM}}\|_{L^2(\Omega)}/\|u_{\varepsilon}^{\text{FEM}}\|_{L^2(\Omega)}$ 
between the TNN solution and the reference solution of 
the problem (\ref{ex_1D_3scale}) for different $\varepsilon$.}\label{tab:1D_3scale_error_eps}
\centering
\begin{tabular}{||c|c|c|c||}
\hline \hline
$\varepsilon$ &$\frac{1}{5}$ & $\frac{1}{10}$ & $\frac{1}{30}$ \\
\hline
Absolute errors & 5.7517e-03 & 3.2778e-03 & 7.5156e-04 \\
\hline
Relative errors & 1.4900e-02 & 8.4442e-03 & 1.9246e-03 \\
\hline \hline
\end{tabular}
\end{table}

\section{Conclusions}
For solving important MPDEs  from practical science and engineers, a type of TNN 
based machine learning method is proposed in this paper. 
Different from other types of neural networks, TNN has the tensor product structure  
which can transform the high dimensional integrations in loss functions to one 
dimensional integrations which can be computed by the classical quadrature schemes 
with high accuracy. This properties lead to high accuracy of the TNN based machine 
learning method for solving the MPDEs. Several numerical examples are 
provided to validate the high accuracy of the proposed methods. 

We believe that the ability of TNN based machine learning method 
can bring more applications in solving more general MPDEs. 
This will be our future work.

\end{document}